\documentclass{amsart}

\usepackage{amsmath, amsthm, amsfonts, amssymb, tikz, mathtools, standalone, reupaper-figs}
\usetikzlibrary{decorations.pathreplacing,angles,quotes}

\DeclarePairedDelimiter{\ceil}{\lceil}{\rceil}

\DeclareMathOperator{\rt}{rt}

\DeclareMathOperator{\urt}{urt}

\newcommand{\abs}[1]{\left\lvert #1 \right\rvert}

\makeatletter 
\newtheorem*{rep@theorem}{\rep@title}
\newcommand{\newreptheorem}[2]{%
\newenvironment{rep#1}[1]{%
 \def\rep@title{#2 \ref{##1}}%
 \begin{rep@theorem}}%
 {\end{rep@theorem}}}
\makeatother

\newtheorem{theorem}{Theorem}
\newreptheorem{theorem}{Theorem}
\newtheorem{lemma}[theorem]{Lemma}

\theoremstyle{definition}
\newtheorem{definition}{Definition}

\begin{document}

\title[Routing convex grid pieces]{Routing by matching on\\ convex pieces of grid graphs}
\author{H.\ Alpert}
\address{Auburn University, 221 Parker Hall, Auburn, AL 36849}
\email{hcalpert@auburn.edu}
\author{R.\ Barnes}
\address{Harvey Mudd College, 320 East Foothill Boulevard, Claremont, CA 91711}
\email{rjbarnes@hmc.edu}
\author{S.\ Bell}
\address{Willamette University, 900 State Street, Salem, OR 97301}
\email{scbell@willamette.edu}
\author{A.\ Mauro}
\address{Stanford University, Building 380, Stanford, CA 94305}
\email{amauro@stanford.edu}
\author{N.\ Nevo}
\address{Colorado College, 14 E.\ Cache La Poudre St., Colorado Springs, CO 80903}
\email{n\_nevo@coloradocollege.edu}
\author{N.\ Tucker}
\address{Juniata College, 1700 Moore Street, Huntingdon, PA 16652}
\email{tuckent18@juniata.edu}
\author{H.\ Yang}
\address{MIT, 77 Massachusetts Avenue, Cambridge, MA 02139}
\email{hannay@mit.edu}
\begin{abstract}
The routing number is a graph invariant introduced by Alon, Chung, and Graham in 1994, and it has been studied for trees and other classes of graphs such as hypercubes.  It gives the minimum number of routing steps needed to sort a set of distinct tokens, placed one on each vertex, where each routing step swaps a set of disjoint pairs of adjacent tokens.  Our main theorem generalizes the known estimate that a rectangular grid graph $R$ with width $w(R)$ and height $h(R)$ satisfies $\rt(R) \in O(w(R)+h(R))$. We show that for the subgraph $P$ of the infinite square lattice enclosed by any convex polygon, we have $\rt(P) \in O(w(P)+h(P))$.
\end{abstract}

\maketitle

\section{Introduction}\label{sec-intro}

Routing number is an invariant of graphs, defined 
by Alon, Chung, and Graham~\cite{Alon94}.  Given a connected graph $G$ on $n$ vertices, we imagine tokens labeled $1$ through $n$ sitting on the vertices of $G$ in some order.  In each routing step, we may select any set of disjoint edges in $G$, and for each edge, swap the tokens on the two vertices of that edge.  Any two of the $n!$ token configurations are connected by some sequence of routing steps; to see this, take a spanning tree of $G$, and move the tokens into position one at a time, starting with the leaves and moving inward.  Thus, we may measure the distance between configurations in terms of routing steps needed.  The \emph{\textbf{routing number}} of $G$, denoted $\rt(G)$, is the maximum, over all pairs of token configurations, of the distance in routing steps between the two configurations.

In our main theorem, the graphs we consider are induced subgraphs of the infinite grid graph, which has vertex set $\mathbb{Z} \times \mathbb{Z}$ and an edge between each pair of vertices with Euclidean distance $1$.  Given a convex polygon $P \subseteq \mathbb{R}^2$, we define the \emph{\textbf{convex grid piece cut out by $P$}} to be the graph $G_P$ with vertices at all lattice points in and on $P$, and edges between pairs of lattice points of distance $1$.  In the remainder of the paper, we use the letter $P$ for both the polygon and the graph, using the notation $\rt(P)$ to mean $\rt(G_P)$.  Although there are some convex polygons $P$ for which the graph $G_P$ is disconnected, the routing number is defined only when $G_P$ is connected.  Note also that when $P$ is translated or rotated, the graph changes, and so the routing number may change slightly.

Our main theorem bounds $\rt(P)$ in terms of the width and height of $P$.  The width $w(P)$ and height $h(P)$ are the maximum absolute differences in $x$-coordinates and in $y$-coordinates, respectively, of any pair of points in $P$.

\begin{theorem}\label{thm-main}
Let $P$ be a connected convex grid piece. Then the routing number of $P$ satisfies the bound $\rt(P) \leq  C (w(P) + h(P))$ for some universal constant $C$. 
\end{theorem}

The reverse inequality is immediate: the diameter of $P$ is within a constant factor of $w(P) + h(P)$, and the routing number of any graph is at least its diameter, because a token may need to travel between two farthest vertices.  Thus we may estimate $\rt(P)$ as $\Theta(w(P) + h(P))$.

One motivation for studying the routing number of convex grid pieces is as a discrete model of configuration spaces of disks.  Given a region $R$ in the plane, such as a convex polygon, the configuration space $\mathrm{Conf}_{n, r}(R)$ as defined in~\cite{Baryshnikov13} is the space of all ways to arrange $n$ disjoint, labeled disks of radius $r$ inside $R$.  If the configuration space is connected, we can define the distance between two configurations to be the amount of time it takes to move between them if the disks can move simultaneously, each with speed at most $1$.  Roughly, the maximum distance between two configurations corresponds to the routing number of the grid piece cut out by $R$; one major difference is that the routing number does not account for what proportion of $R$ is covered by disks, simplifying the problem.

Whereas the routing number of graphs has clear significance in terms of routing information through computer networks, configuration spaces of disks have their own concrete applications.  The $3$-dimensional version of disk configuration spaces is the hard spheres gas model, in which the disks (or spheres) represent individual molecules moving around in a container; see~\cite{Lowen99, Diaconis09} for exposition on the hard spheres model.  If the molecules are densely packed, they can only rattle in place, as in a solid; if there is a lot of space, they can move almost independently, as in a gas, and at intermediate densities the configuration space is somehow like that of a liquid.  Another interpretation of configuration spaces of disks imagines each disk as a robotic car, moving around in an enclosed room such as a warehouse floor.  The geometry and topology of the configuration space constrains what instructions may be used to coordinate the motion of the robots, as in Farber's ``topological complexity''~\cite{Farber08}.

Researchers interested in the robotic car interpretation have made various discrete models of configuration spaces of disks; see, for instance, \cite{Demaine18, Chinta20, Alpert20}.  Typically a discrete result is proved for a rectangular grid, and then the discrete result implies a continuous result about configurations of disks in a rectangular region.  Although restricting attention to a rectangular regions may seem like a minor assumption, the proof structure of the discrete results tends to rely on the rectangular shape.  The reason is that rectangles are self-similar: a rectangular grid is a union of smaller rectangular blocks, with the blocks arranged again in a rectangular grid pattern.  For robotic cars moving in a round disk, for example, these self-similarity properties do not apply.

Thus, the purpose of our theorem is to prove a discrete result for regions that are not necessarily rectangular.  The proof is for convex regions because the claim is not true for arbitrary nonconvex regions; for grid pieces cut out by nonconvex polygons, the bound on routing number is about as bad as for arbitrary trees, which are the hardest to route of all graphs.  Although considering routing number of convex grid pieces is just one possible discrete model for configuration spaces of disks, we hope that the proof method suggests the steps needed to prove such a result for other discrete models as well.

To prove the theorem, we first construct an algorithm for routing tokens on a special class of convex grid pieces, which we call ramp-like polygons.  This class generalizes both rectangles and right triangles, and the recursive algorithm is fairly technical.  Then, we prove that bounds on routing number for some graphs imply bounds on routing number for other graphs: if we can route ramp-like polygons, then we can route polygons cut into two (and then four) ramp-like pieces, and then if we shear these polygons by at most 45 degrees, we can still route the result.  Using these reductions we show the bound for all convex grid pieces.

Section~\ref{sec-augmenting} contains definitions and lemmas needed for the rest of the paper, including the definitions of ramp-like and burger bun polygons.  In Section~\ref{sec-ramp-like}, we prove the routing number bound for the class of ramp-like polygons.  In Section~\ref{sec-burger-bun}, we extend the bound to a more general class which we call burger bun polygons, each of which can be cut into four ramp-like pieces.  Then in Section~\ref{sec-main-proof} we extend the bound to arbitrary convex polygons, using the fact that they can be obtained from burger bun polygons using a shear transformation of at most 45 degrees.

\emph{Acknowledgments.} This research was performed at the MathILy-EST 2020 REU, supported by the National Science Foundation under Award No.\ DMS 1851842.  H.~Alpert was also supported by NSF Award No.\ DMS 1802914.

\section{Preliminaries}\label{sec-augmenting}

In this section, first we give definitions needed for the rest of the paper.  Then we state the known results on routing number that we need.  Finally we prove two lemmas that we use in multiple later sections, showing that adding a small number of vertices to a graph does not increase the routing number by too much.

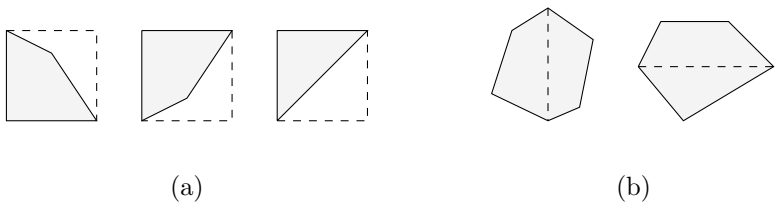
\begin{figure}[h] 
    \centering
\begin{tikzpicture}[scale=.6]
\draw[dashed] (-9,4)--(-7,4)--(-7,2);
\fill[gray!9!] (-9,4)--(-9,2)--(-7,2)--(-8,3.5)--(-9,4);
\draw (-9,4)--(-9,2)--(-7,2)--(-8,3.5)--(-9,4);

\draw[dashed] (-4,4)--(-4,2)--(-6,2);
\fill[gray!9!] (-6,4)--(-4,4)--(-5,2.5)--(-6,2)--(-6,4);
\draw (-6,4)--(-4,4)--(-5,2.5)--(-6,2)--(-6,4);

\draw[dashed] (-3,2)--(-1,2)--(-1,4);
\fill[gray!9!] (-3,2)--(-3,4)--(-1,4)--(-3,2);
\draw (-3,2)--(-3,4)--(-1,4)--(-3,2);

\node[scale=1] at (-5,.5) {(a)};

\fill[gray!9!] (1.75,2.6)--(2.2,4)--(3,4.5)--(4,3.8)--(3.7,2.3)--(3,2)--(1.75,2.6);
\draw (1.75,2.6)--(2.2,4)--(3,4.5)--(4,3.8)--(3.7,2.3)--(3,2)--(1.75,2.6);
\draw[dashed] (3,4.5)--(3,2);
\fill[gray!9!] (5,3.2)--(5.5,4.2)--(7,4.2)--(8,3.2)--(6,2)--(5,3.2);
\draw (5,3.2)--(5.5,4.2)--(7,4.2)--(8,3.2)--(6,2)--(5,3.2);
\draw[dashed] (5,3.2)--(8,3.2);

\node[scale=1] at (4.9,.5) {(b)};
\end{tikzpicture}
    \caption{(a) Ramp-like polygons share at least two sides with their bounding boxes (drawn with dashed lines) and (b) burger-bun polygons either have top and bottom points on the same vertical line, or have leftmost and rightmost points on the same horizontal line.}
    \label{fig:rl_and_bb}
\end{figure}

We define ramp-like polygons and burger bun polygons to be special classes of convex polygons in the plane.  A \emph{\textbf{ramp-like polygon}} is a convex polygon that shares two edges with its bounding box.  That is, there is a rectangle containing our polygon with vertices $(x_1, y_1), (x_1, y_2), (x_2, y_1), (x_2, y_2)$, such that (at least) three of these vertices are vertices of our convex polygon.  A \emph{\textbf{burger bun polygon}} either has top and bottom points on the same vertical line, or has leftmost and rightmost points on the same horizontal line.  That is, either there are two points $(x, y_1)$ and $(x, y_2)$ such that all the $y$-coordinates in the polygon are in the interval $[y_1, y_2]$, or there are two points $(x_1, y)$ and $(x_2, y)$ such that all the $x$-coordinates in the polygon are in the interval $[x_1, x_2]$. Some examples of each of these polygons can be seen in Figure~\ref{fig:rl_and_bb}.

We have defined the routing number $\rt(G)$ of a graph $G$ to be the minimum number of routing steps needed to get from any permutation of labeled tokens on the vertices of $G$ to any other permutation.  Sometimes, instead of having a different label for each token, it helps to consider just two distinct types of tokens, for instance, black tokens and white tokens.  Equivalently, we can consider all of the tokens to be identical, but have some vertices with no tokens on them, so that instead of black tokens and white tokens, we have vertices with tokens and vertices without tokens.  We define the \emph{\textbf{unlabeled routing number}} of $G$, denoted $\urt(G)$, to be the minimum number of routing steps needed to get from any arrangement of black and white tokens, one token per vertex of $G$, to any other arrangement with the same number of black and white tokens.  That is, we take the maximum, over all $k$ and all pairs of arrangements with $k$ black tokens and $\abs{V(G)} - k$ white tokens, of the minimum number of routing steps needed to get from one arrangement to the other.  As before, each routing step consists of selecting a set of disjoint edges of $G$ and swapping the two tokens on the ends of each edge.  Sometimes we refer to routing number as \emph{\textbf{labeled routing number}} to distinguish it from unlabeled routing number.

For reference we state the theorems estimating the routing numbers of paths, trees, and rectangular grids.  The versions that follow are sufficient for our use in this paper, and the proofs can be found in~\cite{Alon94}.

\begin{theorem}[Path bound]\label{thm-path}
For a path $P$ with $n$ vertices, we have $\rt(P) = n$.
\end{theorem}

\begin{theorem}[Tree bound]\label{thm-tree}
For a tree with $n$ vertices, and thus for any connected graph $G$ with $n$ vertices, we have $\rt(G) \leq 3n$.
\end{theorem}

\begin{theorem}[Rectangle bound]\label{thm-rectangle}
For a $p$ by $q$ rectangular grid graph $R_{p, q}$, we have
\[\rt(R_{p, q}) \leq \frac{3}{2}(p+q).\]
\end{theorem}

The two lemmas in the remainder of this section are stated in terms of \emph{\textbf{lattice graphs}}, which we define to be graphs $G$ such that the vertex set is a set of points $(x,y) \in \mathbb{Z} \times \mathbb{Z}$, and there is an edge between two vertices whenever the Euclidean distance between them is exactly $1$.  Our notation sometimes conflates polygons, graphs, and their vertex sets.  When we use set operations on graphs, we typically mean that the operation should be done on the vertex sets, and then we should consider the induced subgraph of $\mathbb{Z} \times \mathbb{Z}$ determined by the resulting set of vertices.  We use the notation $w(G)$ and $h(G)$ for the width and height of a lattice graph $G$, defined similarly to the width and height of a polygon.

In the rest of the paper, we will often trim off, or add on, small or skinny parts of our graphs that do not significantly change the routing numbers.  The following lemma makes that operation precise.

\begin{lemma}\label{lem-hair}
Let $G$ be a connected lattice graph, with its vertices partitioned into sets $K$ (``core'') and $H$ (``hair'').  Suppose that
\begin{enumerate}
\item There are at most $c_1 \cdot (w(G) + h(G))$ vertices in $H$;
\item There is a set of vertices $S \subseteq K$ (``skin'') containing at most $c_2 \cdot (w(G) + h(G))$ vertices, such that the induced subgraph $S \cup H$ is connected; and
\item The routing number of $K$ is at most $c_3 \cdot (w(G) + h(G))$.
\end{enumerate}
Then, we have $\rt(G) \leq (6c_1 + 3c_2 + 2c_3)(w(G) + h(G)).$
\end{lemma}

\begin{proof}

\begin{figure}
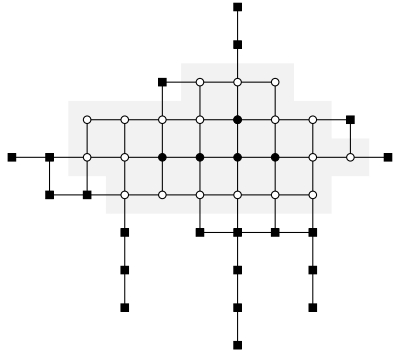

    \centering
    \figcore
    \caption{A partition of a graph $G$ into core $K$ (black and white circular vertices, with gray background) and hair $H$ (black square vertices). The skin $S \subseteq K$ is denoted with white circular vertices.}
    \label{fig:core_hair_partition}
\end{figure}

Figure~\ref{fig:core_hair_partition} shows an example of a core $K$, skin $S \subseteq K$, and the surrounding hair $H$.

First we route within $K$ so that the tokens in $K$ that belong in $H$ have the following property: no token in $K$ that belongs in $H$ has a greater distance to $S$ than a token in $K$ that does not belong in $H$.  That is, for the tokens in $K$ belonging in $H$, we move them to fill $S$ first, then to fill the vertices at distance $1$ from $S$, and so on.  This first phase takes at most $c_3(w(G) + h(G))$ routing steps.

Next, we consider the induced subgraph $S'$ of $G$ containing $H$, $S$, and any other vertices with tokens belonging in $H$.  Because of the previous step and the fact that $H \cup S$ is connected, we know that $S'$ is connected.  It has at most $2\abs{H} + \abs{S} \leq (2c_1 + c_2)(w(G) + h(G))$ vertices.  Thus, we may use at most $3(2c_1 + c_2)(w(G) + h(G))$ routing steps on $S'$ to move all the tokens belonging in $H$ to their home vertices.

Finally, we route within $K$ to move all the tokens belonging in $K$ to their home vertices.  The total number of routing steps is at most $(6c_1 + 3c_2 + 2c_3)(w(G) + h(G))$.
\end{proof}

Sometimes the ``skin'' set $S$ is very easy to describe, but in the final proof we need to be able to find the skin set of an arbitrary convex grid piece.  The following theorem describes how to do so.

\begin{lemma}\label{lem-skin}
Let $P \subseteq \mathbb{R}^2$ be a convex polygon that cuts out the connected lattice graph $K$.  Then there is a connected subgraph $S$ of $K$ containing at most $2(w(K) + h(K))$ vertices, with the following property: if $G$ is a connected lattice graph containing $K$, then $S \cup (G \setminus K)$ is connected.
\end{lemma}

\begin{proof}
We think of $S$ as the circuit enclosing $K$.  To be more precise, we start with the loop $P$, which we may assume is the boundary of the convex hull of $K$.  Then for each edge $E$ of $P$, we modify the loop in the following way.  The two ends of $E$ are lattice points, and we consider the union of grid squares that $E$ passes through.  This union is some centrally symmetric polygon (probably non-convex), and $E$ cuts it into two halves.  

We claim that one half of the boundary of this grid-square polygon is completely contained in the graph $K$.  To see this, suppose to the contrary that another edge $F$ of $P$ also passes through a grid square that $E$ passes through.  Then we can draw a line segment from one point on $E$ to one point on $F$ in the interior of this grid square, and this line segment cuts $P$ into two pieces in a way that separates the two vertices of $E$ but does not intersect the graph $K$.  This contradicts the assumption that $K$ is connected.

Thus, one half of the boundary of the grid-square polygon determined by $E$ is a path in $K$, and in our loop $P$, we may replace $E$ by this path in $K$.  By doing these replacements on all edges of $P$, we obtain a circuit in $K$, and we let $S$ be the set of all vertices and edges in this circuit. $S$ is connected, by following the circuit, as shown in Figure~\ref{fig:circuit_construction}.

\begin{figure}
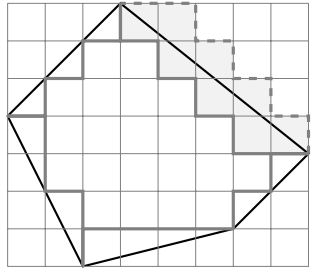

    \centering
    \figcircuit
    \caption{The circuit $S$ is obtained by replacing each edge of $P$ with a path in $K$ which is the boundary of the grid square polygon formed by taking the union of grid squares that $E$ passes through.}
    \label{fig:circuit_construction}
\end{figure}

For each row of vertical edges in $K$, our subgraph $S$ contains only the leftmost and rightmost, and for each column of horizontal edges in $K$, our subgraph $S$ contains only the topmost and bottommost.  Thus the circuit traverses $2h(K)$ vertical edges and $2w(K)$ horizontal edges, and the same number of vertices, so there are at most $2(w(K) + h(K))$ vertices in $S$.

Suppose that $G$ is a connected lattice graph containing $K$, and let $H = G \setminus K$.  Consider any edge $e$ from $H$ to $K$.  We claim that its vertex $v$ in $K$ is in $S$.  In the case where $v$ is on the boundary of the polygon $P$, we know that $v$ is in $S$ because our replacement process to transform from the boundary loop of $P$ to the circuit $S$ does not touch the vertices on the loop.  Otherwise, the edge $e$ crosses from inside $P$ to outside $P$, so it crosses an edge $E$ of $P$.  Then $E$ crosses the two grid squares containing $e$, and so the vertex $v$ on their boundary is part of $S$.
Then, to show that $S \cup H$ is connected, consider any path in $G$ between two vertices of $S \cup H$.  It alternates between sequences of vertices in $H$ and sequences of vertices in $K$, and we have just shown that each $K$ sequence begins and ends with vertices in $S$.  Thus, we may replace each $K$ sequence by an $S$ sequence to get a walk in $S \cup H$ connecting the same two vertices.
\end{proof}

\section{Ramp-like polygons}\label{sec-ramp-like}

The purpose of this section is to prove Theorem~\ref{thm-ramplike}, which bounds the routing number of ramp-like polygons.  Subsection~\ref{subsec:ramps} contains all of the proof except for two big lemmas, which we save for their own subsections: the monotonic configuration theorem (Theorem~\ref{thm:mon}) is proved in Subsection~\ref{subsec:mon}, and the column preparation lemma (Lemma~\ref{lem:rec_rec}) is proved in Subsection~\ref{subsec:cols}.

\subsection{Ramp-like routing overview}\label{subsec:ramps} When routing within ramp-like polygons, the way we use the ramp-like geometry is by defining a slightly more general property of the induced graph, and then using that graph property in the routing.

\begin{definition}
A \textit{\textbf{ramp-like graph}} $R$ is a finite induced subgraph of the infinite lattice graph $\mathbb{Z} \times \mathbb{Z}$ with the following properties:
\begin{enumerate}
    \item Rows are contiguous, and start at $x$-coordinate $0$: if $(x,y) \in V(R)$, then $(i,y) \in V(R)$ for $0 \leq i \leq x$. When we give row numbers, we number the rows by their $y$-coordinates, so that the row numbers increase from bottom to top rather than from top to bottom.
    \item Columns are contiguous, and start at $y$-coordinate $0$: if $(x,y) \in V(R)$, then $(x,i) \in V(R)$ for $0 \leq i \leq y$. When we give column numbers, we number the columns by their $x$-coordinates.
    \item $R$ has discretely convex border: if $n_i$ denotes the greatest $x$-coordinate among all vertices with $y$-coordinate $i$, then for all $i > j \geq c > 0$, we have $n_{i-c} - n_i \geq n_{j-c} - n_j - 1$. 
\end{enumerate}
\end{definition}

\begin{figure}
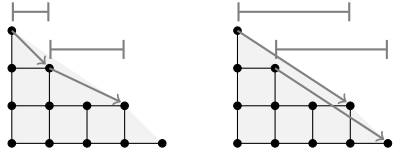

    \centering
    \figdiscconv
    \caption{The discretely convex border property says that the differences between row lengths roughly decrease from bottom to top, failing to decrease by at most 1.}
    \label{fig:discrete-convex}
\end{figure}

The intuition behind the third property is to think of $R$ as being cut out from the first quadrant by the sideways graph $(f(y), y)$ of a function $f$.  If $f$ cuts out a convex shape, then for any constant $c$, the function $f(y-c) - f(y)$ is increasing in $y$.  However, the resulting lattice graph only satisfies the inequality with the term of $-1$ included, as in Figure~\ref{fig:discrete-convex}.

\begin{lemma}
The lattice graph $R$ cut out by any ramp-like polygon $P$ is isomorphic to a ramp-like graph.
\end{lemma}

\begin{proof}
We may assume that $P$ is the convex hull of its enclosed lattice points.  By translating $P$ and rotating by some multiple of a right angle, we may assume that the two sides of $P$ that coincide with sides of its bounding box are along the positive $x$-axis and the positive $y$-axis.  Then the rows and columns are contiguous and start at $0$, so it remains to check the third property, about having discretely convex border.

We consider the rightmost vertices in the rows of $R$ with the $y$-coordinates $i, i -c, j,$ and $j - c$, with $i > j\geq c > 0$. They are $(n_i, i), (n_{i - c}, {i - c}), (n_j, j),$ and $(n_{j - c}, j - c)$, respectively. 

We note that by convexity of the ramp-like polygon $P$, the convex hull of $(n_{i}, i)$, $(n_{j-c}, j- c)$, $(n_i,0)$, $(n_{j - c}, 0)$ is contained in $P$.  We know that $(n_{i-c} + 1, i-c)$ and $(n_j + 1, j)$ are not in $P$, so they also must not be in the convex hull of $(n_{i}, i)$, $(n_{j-c}, j- c)$, $(n_i,0)$, $(n_{j - c}, 0)$. Thus, the lattice points $(n_{i-c} + 1, i-c)$ and $(n_{j} + 1, j)$ must be to the right of the line segment between $(n_{i}, i)$ and $(n_{j-c}, j-c)$.  We let $s$ be the slope of this segment, and compare this slope to the slopes of the segments from $(n_i, i)$ to $(n_{i-c}+1, i-c)$ and from $(n_j+1, j)$ to $(n_{j-c}, j-c)$, as in Figure~\ref{fig:underlying}.

\begin{figure}[h]
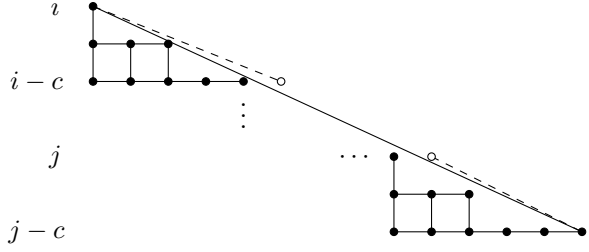

    \centering
    \figunderlying
    \caption{The slope of the solid line, $s$ must be steeper (more negative) than the slope of the upper dashed line, $\frac{-c}{(n_{i-c} + 1) - n_i}$, and shallower (less negative) than the slope of the lower dashed line, $\frac{-c}{n_{j-c} - (n_j + 1)}$.}
    \label{fig:underlying}
\end{figure}

We note that if $s$ is undefined, then each row between $i$ and $j-c$ must contain $n_i = n_{j-c}$ vertices, so the inequality must be true. Otherwise, the slope of the segment from $(n_i, i)$ to $(n_{i-c} +1, i-c)$ is $\frac{-c}{(n_{i-c} + 1) - n_i}$, and the slope of the segment from $(n_j+1), j)$ to $(n_{j-c}, j-c)$ is $\frac{-c}{n_{j-c} - (n_j + 1)}$, and the three negative slopes are ordered as
\[\frac{-c}{n_{j-c} - (n_j + 1)} < s < \frac{-c}{(n_{i-c} + 1) - n_i}.\]
Taking absolute values and comparing the denominators, we have
\[n_{j-c} - (n_j + 1) < (n_{i-c} + 1) - n_i,\]
and since both sides are integers, this is equivalent to our desired inequality
\[n_{i-c} - n_i \geq n_{j-c} - n_j - 1.\]
\end{proof}

The main goal of this section is to prove the following bound on the unlabeled routing number of ramp-like graphs.  This bound implies the corresponding bound for labeled routing number relatively easily.

\begin{theorem}[Unlabeled ramp-like bound]\label{thm-unlabeled}
There is a constant $C$ such that for any ramp-like graph $R$, the unlabeled routing number of $R$ satisfies the bound
\[\urt(R) \leq C \cdot (w(R) + h(R)),\]
where $w(R)$ and $h(R)$ denote the width and height of $R$ as a lattice graph.
\end{theorem}

Given any $k$, we define the \emph{\textbf{row-major order}} configuration of $k$ black tokens and $\abs{V(R)} - k$ white tokens on the graph $R$ as follows: we order the vertices in $R$ by row from top to bottom, and within each row from left to right, and we take the configuration in which all of the black tokens appear before all of the white tokens.  To prove the unlabeled ramp-like bound, we start with an arbitrary configuration of black and white tokens on $R$, and describe how to route from this configuration to row-major order.  

The process for arbitrary ramp-like graphs is considerably more complicated than it is for rectangular grids.  On a rectangular grid, given a configuration of black and white tokens, we can move them within their rows to get the right number of each color into each column, and then move them within their columns to achieve row-major order.  An arbitrary ramp-like graph may be much narrower at the top than at the bottom, so our starting configuration could be a few wide rows of black tokens along the bottom, which we want to move to form several narrow rows at the top.  In this case we would need to alternate between horizontal and vertical motion several times to move between the configurations.  Whatever the starting configuration is, our first phase is to route to what we call a monotonic configuration.

\begin{figure}[h]
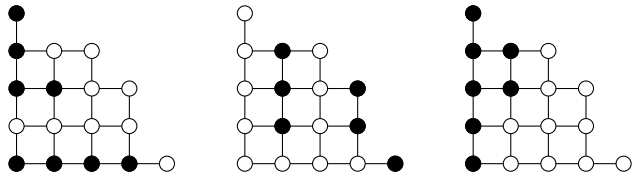

    \centering
    \figlauamon
    \caption{(a) A left-aligned configuration of tokens, (b) an up-aligned configuration of tokens, and (c) a monotonic configuration of tokens.}
    \label{fig:la_ua_monotonic}
\end{figure}

\begin{definition}
We say a configuration of tokens on a ramp-like graph $R$ is \emph{\textbf{left-aligned}} if for every black token on some $(x,y) \in V(R)$, there is also a black token on $(x - 1, y) \in V(R)$, or $(x - 1, y) \notin V(R)$. Similarly, we say a configuration of tokens is \emph{\textbf{up-aligned}} if for every black token on some $(x,y) \in V(R)$, either there is also a black token on $(x, y + 1) \in V(R)$, or $(x, y + 1) \notin V(R)$. Additionally, we say a given black token is left-aligned or up-aligned if it satisfies the corresponding condition stated above.

We say a configuration is \emph{\textbf{monotonic}} if it is both left-aligned and up-aligned. Examples of configurations with each of these properties can be seen in Figure~\ref{fig:la_ua_monotonic}.
\end{definition}

We note that routing all black tokens such that they are as far up or left as possible within their column or row results in an up-aligned or left-aligned configuration, respectively. 
We abbreviate the process of routing all black tokens as far up as possible (within their column) with the phrase ``pushing up,'' and routing all black tokens as far left as possible (within their row) with ``pushing left.''

The following theorem describes a process sufficient for moving the tokens to a monotonic configuration.

\begin{theorem}[Monotonic configuration]\label{thm:mon}
For any configuration $X_0$ of black and white tokens on a ramp-like graph, after pushing the black tokens up, then left, then up, then left, the new configuration of tokens is monotonic. 
\end{theorem}

Proving the monotonic configuration theorem (Theorem~\ref{thm:mon}) is the most technical aspect of routing within ramp-like polygons, and the proof appears in Subsection~\ref{subsec:mon}.  After routing to a monotonic configuration, we are in a better position to move the tokens to the correct columns, after which all that remains is to push up all the black tokens to get to row major order.  Our proof that we can move the tokens to the correct columns sufficiently quickly is by induction on the number of rows in our graph.  To do this, we want to partition the graph into a top slice (denoted by $R_{m_1 - 1}$ in the next lemma) and a bottom slice and route the two slices separately; that is, we want a way to move the tokens to the correct columns without crossing  between the two slices.  The next lemma states that this division into top and bottom is possible.  For the lemma statement, we denote the number of vertices in the $y=i$ row by $1+n_i$, so that as before the $x$-coordinates of those vertices range from $0$ to $n_i$.

\begin{lemma}[Column preparation]\label{lem:rec_rec}
Let $R$ be a ramp-like graph with $m$ rows (that is, from $y=0$ through $y=m-1$), and let $R_i$ denote the induced subgraph of $R$ consisting of the top $i$ rows (that is, from $y=m-i$ through $y = m-1$).  Let $X$ be a monotonic configuration on $R$ with $t$ black tokens, and let $m_1$ be such that
\[ \sum^{m_1 - 1}_{i = 1} (1+n_{m-i}) < t \leq \sum^{m_1}_{i = 1} (1+n_{m-i}),\]
meaning that $t$ black tokens can fit on the vertices of $R_{m_1}$, but are not able to fit on the vertices of $R_{m_1 - 1}$.  Then, there exists a configuration $Y$ on $R$ such that 
\begin{enumerate}
    \item $R_{m_1 - 1}$ contains the same number of tokens of each color in $Y$ as in $X$, and
    \item In $Y$, each column of $R$ contains the same number of tokens of each color as there are in row major order.
\end{enumerate}
\end{lemma}

The proof of the column preparation lemma (Lemma~\ref{lem:rec_rec}) is also fairly technical, and it appears in Subsection~\ref{subsec:cols}.  Outside of the proof of the lemma, we do not need to remember the definition of $m_1$; rather, what is important is the conclusion of the lemma, which makes $R_{m_1 - 1}$ the top slice and the remainder of $R$ the bottom slice.  It would be more intuitive to choose the top slice to be $R_{m_1}$, but then in the case where $R_{m_1}$ is the whole graph $R$, we would not be able to apply the inductive hypothesis to it.  Thus, we choose the top slice to be $R_{m_1 - 1}$.  Assuming the monotonic configuration theorem (Theorem~\ref{thm:mon}) and the column preparation lemma (Lemma~\ref{lem:rec_rec}), we can finish proving the upper bound on unlabeled routing number of ramp-like graphs.

\begin{proof}[Proof of unlabeled ramp-like bound (Theorem~\ref{thm-unlabeled})]
Let $m$ be the number of rows in our ramp-like graph $R$, and let $n$ be the number of columns.  If $m$ or $n$ is $1$, then $R$ is a path, which we already know how to route by the path bound (Theorem~\ref{thm-path}).  Thus, we may assume that $w(R)$ and $h(R)$ are both at least $1$.  In this case, we have $n = w(R) + 1 \leq 2w(R)$ and $m = h(R) + 1 \leq 2h(R)$, so it suffices to find a constant $C$ such that we can move an arbitrary configuration of black and white tokens to row-major order in at most $C(m+n)$ routing steps.

Let $C_0$ be the constant from the rectangle bound (Theorem~\ref{thm-rectangle}) such that a rectangular grid with $p$ rows and $q$ columns has unlabeled routing number at most $C_0(p + q)$.

First we show by induction on $m$ that the process in the column preparation lemma (Lemma~\ref{lem:rec_rec}) of moving from configuration $X$ to configuration $Y$ can be accomplished in at most $C_0(m + n)$ moves.  In the base case $m = 1$, our ramp-like graph is a path of length $n$, which is already rectangular. For $m > 1$, we partition the ramp-like graph $R$ into three parts: the upper ramp-like shape $R_{m_1 - 1}$, the rectangular piece containing the row below $R_{m_1 - 1}$ and all vertices directly below this row, and the ramp-like graph consisting of columns to the right of the rectangular piece.  

\begin{figure}
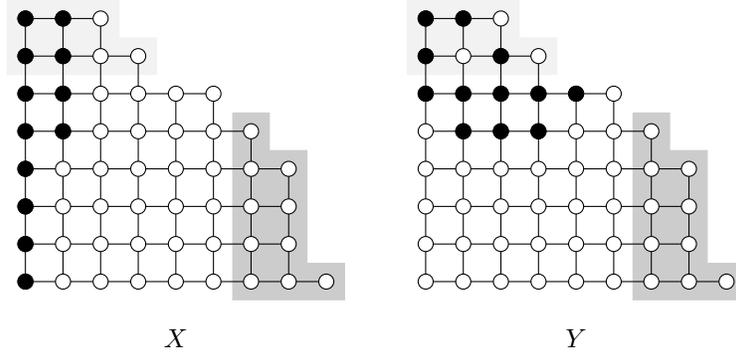

\centering 
\figpartitionrl
\caption{$R$ is partitioned into $R_{m_1 - 1}$ (light gray), a rectangle, and a subgraph to the right of the rectangle (dark gray).  $R_{m_1 - 1}$ is routed using the inductive hypothesis, and in parallel, the rectangular piece is routed to produce a configuration, which, when pushed up, is in row major order.  The subgraph to the right of the rectangle contains no black tokens either in the monotonic configuration $X$ or in row-major order.}
\label{fig:partition_rl}
\end{figure}

To move from configuration $X$ to configuration $Y$, we apply the inductive hypothesis to $R_{m_1 - 1}$ while simultaneously routing within the rectangular piece; the graph to the right of the rectangular piece has no tokens in it, so it does not need any routing.  Because $R_{m_1 - 1}$ has strictly fewer rows than the original graph $R$, the inductive hypothesis applies.  The rectangular piece has at most $m$ rows and at most $n$ columns, so it requires at most $C_0(m + n)$ routing steps to move from $X$ to $Y$, and by the inductive hypothesis, $R_{m_1 - 1}$ also requires at most $C_0(m+n)$ routing steps.  Performing the steps simultaneously completes the induction.

Given an arbitrary token configuration $X_0$ on the ramp-like graph $R$, we start by pushing all the tokens up, then left, then up, then left, which by the monotonic configuration theorem (Theorem~\ref{thm:mon}) gives a monotonic configuration.  We apply the column preparation lemma (Lemma~\ref{lem:rec_rec}) to the resulting monotonic configuration, and then we push all the tokens up to get to row-major order.  By the path bound (Theorem~\ref{thm-path}) it takes at most $m$ routing steps to push up and at most $n$ routing steps to push left, so the total number of routing steps to get  to row-major order is at most $3m + 2n + C_0(m+n) \leq (3 + C_0)(m+n)$.

Given any two configurations of the same set of white and black tokens on $R$, we can route from one to the other by routing the first into row-major order, and then routing from row-major order to the second.  Thus, if we choose $C = 4(3+C_0)$, we can route between the two configurations in $2(3+C_0)(m+n) \leq C(w(R) + h(R))$ routing steps.
\end{proof}

We can use the bound on unlabeled routing number of ramp-like graphs to give a bound on labeled routing number, by dividing the graph into four quadrants and applying recursion.

\begin{theorem}[Labeled ramp-like bound]\label{thm-ramplike}
There is a constant $C$ such that for any ramp-like polygon $P$, the routing number of $P$ satisfies the bound
\[\rt(P)\leq C\cdot (w(P)+h(P)).\]
\end{theorem}

\begin{proof}
The idea of the proof is to divide $P$ into quadrants, each with width at most $\frac{1}{2}w(P)$ and height at most $\frac{1}{2}h(P)$.  Using the unlabeled ramp-like bound (Theorem~\ref{thm-unlabeled}), we can move each token into the correct quadrant.  Then, each quadrant is ramp-like, so we can apply recursion to route within all four quadrants simultaneously, as shown in Figure~\ref{fig:recursive_route}.

\begin{figure}[h]
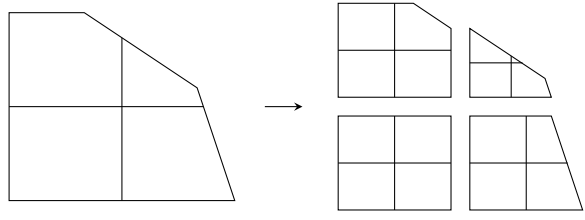

\centering 
\figrecursiveroute
\caption{$P$ is divided into four quadrants, and all tokens are routed into the correct quadrant using unlabeled routing twice. Then, this process is repeated recursively for each quadrant, simultaneously.}
\label{fig:recursive_route}
\end{figure}

We may assume that the vertical and horizontal sides that $P$ shares with its bounding box have rational $x$-coordinate and $y$-coordinate, respectively, and that the width and height of $P$ are both irrational.  These assumptions guarantee that when we cut $P$ in half, the cut does not go through any lattice points.  If $P$ does not have these properties already, we can make $P$ very slightly bigger so that it does, which increases the right-hand side of the desired inequality very slightly; taking the limit of a shrinking sequence of approximations gives the desired inequality.

Let $C_0$ be the constant for unlabeled routing from the unlabeled ramp-like bound (Theorem~\ref{thm-unlabeled}).  We divide $P$ into halves with a vertical line bisecting the width, and then into quadrants with a horizontal line bisecting the height.  One of the quadrants may be empty.  It takes at most $C_0(w(P) + h(P))$ routing steps to move the tokens so that those that belong in the left half go to the left half, and those that belong to the right half go to the right half.  Then each half is ramp-like and has width $\frac{1}{2}w(P)$ and height $h(P)$, so it takes at most $C_0\left(\frac{1}{2}w(P) + h(P)\right)$ additional routing steps to move the tokens into the quadrants where they belong, working with both halves simultaneously.

We select $C = 4C_0$, and prove the theorem by induction on $\ceil{w(P) + h(P)}$.  If $\ceil{w(P) + h(P)} = 1$, then $P$ has only one vertex, so the routing number is $0$, which is certainly at most $C \cdot (w(P) + h(P))$.  Otherwise, we have $\ceil{\frac{1}{2}w(P) + \frac{1}{2}h(P)} < \ceil{w(P) + h(P)}$, so we may apply the inductive hypothesis to find that the routing number of each quadrant of $P$ is at most $C \cdot \left(\frac{1}{2}w(P) + \frac{1}{2}h(P)\right)$.  Then the total number of steps to route an arbitrary configuration of tokens on $P$ to a home configuration is at most
\begin{align*}
C_0(w(P) + h(P)) &+ C_0\left(\frac{1}{2}w(P) + h(P)\right) + C\left(\frac{1}{2}w(P) + \frac{1}{2}h(P)\right)\\
&\leq \left(2C_0 + \frac{1}{2}C\right)\cdot(w(P) + h(P))\\
&= C\cdot(w(P) + h(P)).
\end{align*}
\end{proof}

\subsection{Moving to a monotonic configuration}\label{subsec:mon}

In this subsection we prove the monotonic configuration theorem (Theorem~\ref{thm:mon}).  To understand the strategy, we observe that if our ramp-like graph were a rectangular grid, then pushing the black tokens up and then left would already give a monotonic configuration.  This is because after pushing up from an arbitrary configuration, every black token below the top row has another black token directly above it, so the number of tokens in each row is non-increasing as we consider the rows from top to bottom.  However, for an arbitrary ramp-like graph, after pushing up from an arbitrary configuration, a row might have more black tokens than the row above it, because the lower row might have black tokens in columns to the right of all columns in the upper row.

The proof of the theorem is based on Lemmas~\ref{lem:exist-mount} and~\ref{lem:comp-mount}, which together show that after pushing up, left, and up on a ramp-like graph, the result is similar to what we would get from simply pushing up on a rectangular grid.  Namely, if we consider the rightmost column that contains black tokens, then we show that all of the rows that are too short to extend to that column are completely full of black tokens.  The configuration below these rows looks like an up-aligned configuration on a rectangular grid, so pushing left one more time results in a monotonic configuration.

Lemma~\ref{lem:exist-mount} shows how the shape of our ramp-like graph affects the possible numbers of black tokens per row after our first step of pushing up, and thus after our second step of pushing left as well.  Then Lemma~\ref{lem:comp-mount} describes the result of our third step of pushing up.  For both lemmas, we use the following notation.
Given a ramp-like graph $R$, we let $\#(R, i)$ denote the number of vertices in the row of $R$ with $y$-coordinate $i$.  Given a configuration $X$ on $R$, we let $\#(X, i)$ denote the number of black tokens in $X$ with $y$-coordinate $i$.

\begin{lemma}\label{lem:exist-mount}
Let $X_1$ be an up-aligned configuration on a ramp-like graph $R$.  Then for all rows $b$ and $b+c > b$ of $R$, we have
\[\#(X_1, b) - \#(X_1, b+c) \leq \#(R, b) - \#(R, b+c).\]
\end{lemma}

\begin{proof}
Let $d$ be the number of black tokens in row $b$ of $X_1$ that are to the right of column $n_{b+c}$, the rightmost column of row $b+c$.  Using the fact that $X_1$ is up-aligned, we have
\[\#(X_1, b) - \#(X_1, b+c) \leq d,\]
because every token in row $b$ has a token directly above it in row $b+c$, except for those in the $d$ columns to the right of column $n_{b+c}$.  Because the total number of columns in row $b$ to the right of $n_{b+c}$ is $\#(R, b) - \#(R, b+c)$, we also have
\[d \leq \#(R, b) - \#(R, b+c).\]
Together, these inequalities give the desired inequality
\[\#(X_1, b) - \#(X_1, b+c) \leq \#(R, b) - \#(R, b+c).\]
\end{proof}

Geometrically, the lemma says that after we push up to form $X_1$ and then left to form a configuration $X_2$, the right boundary of the cluster of black tokens is steeper than the right boundary of the graph $R$.  The next lemma starts with this configuration $X_2$ that results from pushing up and left, and describes what happens after pushing up again.

\begin{lemma}\label{lem:comp-mount}
Let $\tau$ be a black token in a left-aligned configuration $X_2$ on a ramp-like graph $R$.  Let $(x_2, b)$ be the coordinates of $\tau$, and let $r > b$ be a row number such that $n_r < x_2$, if such a row exists.  Suppose we know that for all $c > 0$ such that $b+c$ is a row of $R$, we have
\[\#(X_2, b) - \#(X_2, b+c) \leq \#(R, b) - \#(R, b+c).\]
Then when we push up to reach another configuration $X_3$, all rows $r$ and above will have only black tokens.
\end{lemma}

\begin{proof}
We imagine translating all of the tokens upward in the lattice so that row $b$ moves up to row $r$ and some tokens may occupy lattice points that are not in the graph $R$.  If all rows $r$ and above are covered by black tokens in this arrangement, they also have only black tokens in $X_3$.  Thus, it suffices to show that for all $c \geq 0$, the number of black tokens in row $b+c$ of $X_2$ is at least the number of vertices in row $r+c$.

To show this, we combine the inequalities in the hypothesis with an inequality resulting from the discretely convex border property of $R$: because $r > b$, we have
\[n_r - n_{r+c} \geq n_b - n_{b+c} - 1,\]
or equivalently,
\[\#(R, r) - \#(R, b+c) + 1 \geq \#(R, b) - \#(R, b+c).\]
Putting the inequalities together, we have
\begin{align*}
    \#(X_2, b+c) &= \#(X_2, b) - [\#(X_2, b) - \#(X_2, b+c)]\\
    &> \#(R, r) - [\#(R, b) - \#(R, b+c)]\\
    &\geq \#(R, r) - [\#(R, r) - \#(R, r+c) + 1]\\
    &= \#(R, r+c) - 1,
\end{align*}
and so because all of the quantities are integers, we obtain our desired inequality $\#(X_2, b+c) \geq \#(R, r+c)$.
\end{proof}

Having described the configuration that results from pushing up, left, and up, we are ready to prove that pushing this configuration left results in a monotonic configuration.

\begin{proof}[Proof of monotonic configuration theorem (Theorem~\ref{thm:mon})]

We label the sequence of configurations as follows: let
\begin{itemize}
    \item the starting configuration be $X_0$,
    \item the configuration after pushing up be $X_1$,
    \item the configuration after pushing left be $X_2$,
    \item the configuration after pushing up a second time be $X_3$, and
    \item the configuration after pushing left a second time be $X_4$, which is also the final configuration.
\end{itemize}

Note that $X_4$ is left-aligned, so in order to show $X_4$ is monotonic, we just need to show it is up-aligned. We will do this by showing that all black tokens in $X_4$ are up-aligned.

We consider an arbitrary black token $\tau$ in $X_0$. Let the vertex which $\tau$ is on in $X_2$ be $(x_2, b)$. By Lemma~\ref{lem:exist-mount}, we have the inequality
\[\#(X_1, b) - \#(X_1, b+c) \leq \#(R, b) - \#(R, b+c)\]
for all rows $b < b+c$ of $R$, and because every token is in the same row in $X_2$ as in $X_1$, the same inequality is true of $X_2$.
Then by Lemma~\ref{lem:comp-mount}, we have that when we push up to get to $X_3$, rows $r$ and above will be all black.

Let $s$ be the row such that $\tau$ is in row $s-1$ in $X_3$.  Because $\tau$ is in the same column $x_2$ in $X_2$ and $X_3$, and this column is to the right of $n_r$, we have $s - 1 < r$, or in other words $s \leq r$.  If $s = r$, then $\tau$ is up-aligned in $X_4$ because the row above $\tau$ is row $r$, and we have shown that rows $r$ and above are all black in $X_3$.

If $s < r$, then because of how $r$ is defined we have $x_2 \leq n_s$.  Thus, in $X_3$ (which is up-aligned), for $\tau$ and every black token to the left of it in row $s-1$, there is a corresponding black token immediately above, in row $s$.  When we push left to get $X_4$, there are at least as many black tokens in row $s$ as there are black tokens in row $s-1$ to the left of and including $\tau$, so $\tau$ is up-aligned.

Thus, for any black token $\tau$ in $X_0$, in $X_4$, there is a black token (or no vertex) above it and a black token (or no vertex) to the left of it.
\end{proof}

\subsection{Distributing tokens among columns}\label{subsec:cols}

Proving the column preparation lemma (Lemma~\ref{lem:rec_rec}) is the last piece needed to complete the proof of the unlabeled ramp-like bound (Theorem~\ref{thm-unlabeled}) and thus the labeled ramp-like bound (Theorem~\ref{thm-ramplike}).  The goal is to get the right number of tokens of each color into each column, without moving tokens between the top slice $R_{m_1 - 1}$ and the bottom slice $R \setminus R_{m_1 - 1}$.  For the top slice, we choose to move the tokens to row-major order; this determines how many black tokens we want in each column of the bottom slice.  The only thing that could potentially go wrong is if we have somehow assigned more black tokens to a column of the bottom slice than its number of vertices.  We show this does not happen, roughly because the black tokens in the bottom slice are more evenly spaced, among at least as many columns, in our target configuration $Y$ as in our starting configuration $X$.

\begin{proof}[Proof of column preparation lemma (Lemma~\ref{lem:rec_rec})]
Let $t_1$ be the number of black tokens in $R_{m_1 - 1}$ in configuration $X$, and let $m_0 \leq m_1 - 1$ be the number such that 
\[ \sum^{m_0 - 1}_{i = 1} (1+n_{m-i}) < t_1 \leq \sum^{m_0}_{i = 1} (1+n_{m-i}) ,\]
meaning that $t_1$ black tokens can fit on the vertices of $R_{m_0}$, but are not able to fit on the vertices of $R_{m_0 - 1}$.

On any ramp-like graph, we denote the configuration of $t$ black tokens (and the remainder white tokens) in row-major order by $RM(t)$.  We set $Y$ to be equal to $RM(t_1)$ on $R_{m_1-1}$.  Let $Z$ be the configuration of $t-t_1$ black tokens, one at every vertex where $RM(t)$ has a black token but $RM(t_1)$ does not, and let $z_j$ be the number of black tokens in $Z$ that are in the $x=j$ column of $R$.  We note that $z_j = 0$ for $j > n_{m_1}$.  To prove the lemma, we need to show that $z_j$ black tokens can fit into column $j$ of $R \setminus R_{m_1-1}$; that is, $z_j \leq m - m_1 + 1$.

First we address the case where $X$ contains a black token in a column strictly to the right of the subgraph $R_{m_0}$.  We claim that in this case, $z_j \leq 1 \leq m - m_1 + 1$ for all $j$, so we are done.  Because $X$ is monotonic, if $X$ contains a black token to the right of $R_{m_0}$, then $R_{m_0}$ must be entirely full of black tokens, and the row below it must also contain black tokens.  Thus, $R_{m_0 + 1}$ contains more than $t_1$ black tokens.  Because $R_{m_1 - 1}$ contains only $t_1$ black tokens, this then implies that $R_{m_1 - 1}$ must be smaller than $R_{m_0 + 1}$.  Because we always have $m_0 \leq m_1 - 1$, we conclude that in this case we have $m_1 - 1 = m_0$, with $R_{m_1 - 1}$ entirely full of black tokens.  The definition of $m_1$ implies that the black tokens not in $R_{m_1 - 1}$ all fit into the row just below $R_{m_1 - 1}$, so $z_j \leq 1$.

Thus, we may assume that we are in the case where all black tokens in $X$ are in columns $0$ through $n_{m_0}$.  In this case, the idea of the proof is that if we were to distribute the tokens in the bottom slice $R \setminus R_{m_1 - 1}$ as evenly as possible among the columns $0$ through $n_{m_0}$, the column with the most black tokens would have at least as many as in $Z$, because $Z$ may use the columns to the right of $n_{m_0}$ as well.

More precisely, the portion of $X$ in $R \setminus R_{m_1-1}$ has $t - t_1$ black tokens, all of which are in columns $0$ through $n_{m_0}$, so we have
\begin{equation}\label{eqn}
t - t_1 \leq (1+n_{m_0})(m-m_1+1).
\end{equation}
Also, because $Z$ has $t - t_1$ black tokens in total, we have
\[t - t_1 = \sum_{0 \leq j \leq m_1} z_j.\]
Let $\delta$ be the greatest $x$-coordinate of the black tokens in the $y=m-m_1$ row of $RM(t)$ (the $m_1$st row from the top, and the last not-all-white row), and let $\delta_1$ be the greatest $x$-coordinate of the black tokens in the $y=m-m_0$ row of $RM(t_1)$ (the $m_0$th row from the top, and the last not-all-white row).  We now estimate $z_j$ in the three possible cases for how $\delta_1$ and $\delta$ compare, depicted in Figure~\ref{fig-cases}.

\begin{figure}
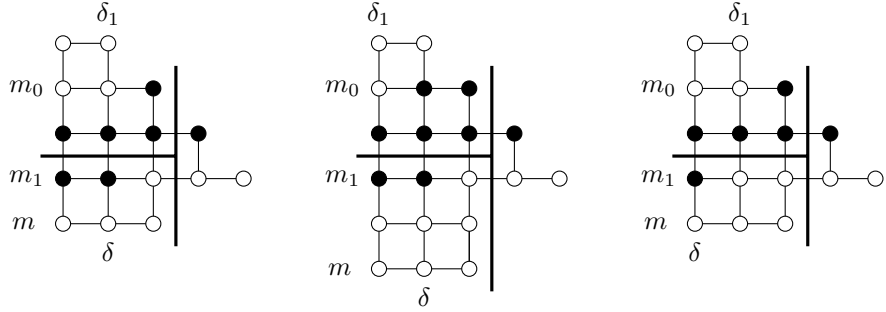

\begin{center}
\figroom     
\end{center}
\caption{In every case, the black tokens of $Z = RM(t) \setminus RM(t_1)$ have enough vertical space to slide down below subgraph $R_{m_1 - 1}$, because the total area of columns $0$ through $n_{m_0}$ in $R \setminus R_{m_1 - 1}$ is large enough for the black tokens in columns $0$ through $n_{m_0}$ of $Z$.}\label{fig-cases}
\end{figure}

If $\delta_1 = \delta$, or if $n_{m_0} = \delta_1 < \delta$, then every column $j$ with $j \leq n_{m_0}$ has $z_j = m_1 - m_0$.  For $j > n_{m_0}$, because row $m-m_0$ has no vertices in column $j$, we have $z_j \leq m_1 - m_0$.  The total number $t - t_1$ of black tokens in $Z$ is at least the number in columns $0$ through $n_{m_0}$, so we have $t - t_1 \geq (1+n_{m_0})(m_1 - m_0)$.  Combining this inequality with inequality~(\ref{eqn}), we have
\[(1+n_{m_0})(m_1 - m_0) \leq t - t_1 \leq (1+n_{m_0})(m - m_1+1),\]
so $m_1 - m_0 \leq m - m_1+1$, and so $z_j \leq m - m_1+1$ for all $j$.

If $\delta_1 < \delta$ and $\delta_1 < n_{m_0}$, then for $j \leq \delta_1$ and for $\delta < j \leq n_{m_0}$ we have $z_j = m_1 - m_0$, and for $\delta_1 < j \leq \min(\delta, n_{m_0})$ we have $z_j = 1 + m_1 - m_0$.  For $j > n_{m_0}$, we have $z_j \leq 1 + m_1 - m_0$ as well.  Thus we have
\[(1+n_{m_0})(m_1 - m_0) < t - t_1 \leq (1+n_{m_0})(m - m_1+1),\]
which implies that $m_1 - m_0 < m - m_1 + 1$, and so $z_j \leq m - m_1+1$ for all $j$.

If, finally, $\delta_1 > \delta$, then for $j \leq \delta$ and for $\delta_1 < j \leq n_{m_0}$ we have $z_j = m_1 - m_0$, and for $\delta < j \leq \delta_1$ we have $z_j = -1 + m_1 - m_0$.  For $j > n_{m_0}$, we have $z_j \leq m_1 - m_0$ as well.  Thus we have
\[(1+n_{m_0})(-1 + m_1 - m_0) < t - t_1 \leq (1+n_{m_0})(m - m_1+1),\]
which implies that $-1 + m_1 - m_0 < m - m_1 + 1$, and so $z_j \leq m - m_1+1$ for all $j$.

Thus, in every case, $R \setminus R_{m_1-1}$ has enough space to fit the same number of black tokens in each column as $Z$.  We set $Y$ to be any configuration that in $R_{m_1 - 1}$ has the same number of black tokens in each column as $RM(t_1)$, and that in $R \setminus R_{m_1 - 1}$ has $z_j$ black tokens in column $j$ for each $j$.
\end{proof}

\section{Burger bun polygons}\label{sec-burger-bun}

In this section we prove Theorem~\ref{thm-burger-bun}, the bound on routing number of burger bun polygons.  Our strategy is to divide the burger bun polygon in half, then to divide each half into two ramp-like pieces.  We know that we can route a single ramp-like piece from the ramp-like bound (Theorem~\ref{thm-ramplike}), so our first step is to use this to prove that we can route a pair of ramp-like pieces and thus a half of a burger bun polygon.  Then, using a similar argument, we show that this implies that we can route a whole burger bun polygon.

We have defined a ramp-like polygon to be a convex polygon such that two of its edges coincide with edges of its bounding box; we refer to each of these edges as a \textit{\textbf{spine}} of the ramp-like polygon.  If two otherwise disjoint ramp-like polygons have a common spine, then their union is also a convex polygon.  Similarly, every burger bun polygon is divided in two by a spine.  Specifically, if a burger bun has two points of maximum and minimum $y$-coordinate with equal $x$-coordinate, then we refer to the segment between those two points as its (vertical) spine, and if it has two points of maximum and minimum $x$-coordinate with equal $y$-coordinate, then we refer to the segment between those two points as its (horizontal) spine.

Our strategy for this section is as follows.  In our situation, we have two polygons with a common spine, and we may assume that we know how to route within each of the two polygons.  Given an arbitrary configuration of labeled tokens on the union of the polygons, we want to route those tokens to their home positions.  It suffices to get each token into the half where it belongs, because then we can route within the halves separately to get each token to its home vertex.  Thus, we have an unlabeled routing problem, thinking of the tokens belonging in the first half as black, and the tokens belonging in the second half as white.

In the special case where the two polygons, and their corresponding graphs, are mirror images across the spine, we can route as follows: first we route within the second half so that the configuration is a color-reversed mirror image of the first half.  This is possible because when the two halves have the same number of vertices, the number of black tokens in the second half is equal to the number of white tokens in the first half.  Once the two halves are color-reversed mirror images, each row (if the spine is vertical) has the same number of white tokens as black tokens, so we may route all rows simultaneously to get each token into the half where it belongs.

\begin{figure}
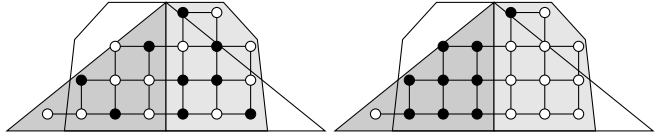

    \centering
    \figmirror
    \caption{We think of the tokens belonging in the left half as black, and in the right half as white.  We can make the right subpolygon a color-reversed mirror image of the left, and then we can route within rows to move the subpolygon tokens to their correct halves.}
    \label{fig:mirror}
\end{figure}

In the general case, where the two polygons are two ramp-like polygons or two halves of a burger bun polygon, the polygons may not be mirror images.  Instead, we find smaller polygons inside them that are mirror images and still contain a significant fraction of the vertices.  Then we can solve the unlabeled routing problem by repeatedly applying the mirror-image technique to these smaller polygons; Figure~\ref{fig:mirror} illustrates this strategy.  The following theorem finds those smaller polygons that are mirror images, for the case where the two original polygons are ramp-like polygons.  We consider both the case where the ramp-like polygons together form half of a burger bun (that is, their non-shared spines are collinear) and the case where they do not, because this latter case turns out to be useful in the next part of the proof, where the two polygons are halves of a burger bun.

\begin{theorem}[Intersection magnitude]\label{thm-int-mag}
Let $P_1$ and $P_2$ be two ramp-like polygons with common vertical spine $E$, such that the widths of $P_1$ and $P_2$ and the length of $E$ are all at least $41$.  Then there exist subgraphs $G_1$ and $G_2$ of $P_1$ and $P_2$, respectively, both disjoint from $E$, such that $G_2$ is a reflection of $G_1$ over some vertical line, and the equal number of vertices in $G_1$ or $G_2$ is at least $\frac{1}{20} \min\{|P_1|, |P_2|\}$, where the absolute value bars denote the number of lattice points in the interior and boundary of each polygon.
\end{theorem}

The proof of this theorem is based on two lemmas: the spine alignment lemma (Lemma~\ref{lem-spine}), and the triangle trimming lemma (Lemma~\ref{lem-trim}).  Roughly, the idea is that to find the subgraphs $G_1$ and $G_2$, we should reflect $P_1$ over the shared spine and intersect it with $P_2$ to find $G_2$, or reflect $P_2$ over the shared spine and intersect it with $P_1$ to find $G_1$.  The spine alignment lemma (Lemma~\ref{lem-spine}) accounts for the fact that the shared spine might not be at an integer or half-integer coordinate, so reflecting across it might not take lattice points to lattice points.  Then the triangle trimming lemma (Lemma~\ref{lem-trim}) starts from a quick estimate of the area of the polygon intersection, and produces an estimate of the number of lattice points inside that polygon.

\begin{lemma}[Spine alignment]\label{lem-spine}
Let $P$ be a burger bun polygon with vertical spine $E$ dividing $P$ into left side $P_1$  and right side $P_2$.  Suppose that the $x$-coordinate of $E$ is not an integer.  Then there is another burger bun $P'$ with vertical spine $E'$, dividing $P'$ into left side $P'_1$ and right side $P'_2$, with the following properties:
\begin{itemize}
\item The $x$-coordinate of $E'$ is an integer.
\item $P_1$ and $P'_1$ contain the same lattice points.
\item The set of lattice points inside $P'_2$ is obtained by translating the set of lattice points inside $P_2$ one unit to the right.
\end{itemize}
\end{lemma}

\begin{figure}[h]
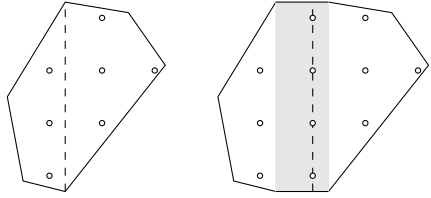

    \centering
    \figspine
    \caption{If we cut a burger bun polygon along its spine and insert a rectangle of width $1$ instead, the result is a burger bun polygon with one additional column and a spine along that column.}
    \label{fig:spine-align}
\end{figure}

\begin{proof}
Figure~\ref{fig:spine-align} shows the relationship between the polygons.  Let $P_1$ and $P_2$ be arbitrary convex polygons sharing a vertical edge $E$ with a non-integer $x$-coordinate $x_0$. Construct a new vertical line segment $E'$ with the same length as $E$ and integer $x$-coordinate $\lceil x_0 \rceil$. Translate all polygon vertices of $P_2$ one unit to the right, and take the convex hull of these vertices with the endpoints of $E'$ to form congruent polygon $P_2'$. Similarly, take the convex hull of all polygon vertices of $P_1$ with the endpoints of $E'$ to form congruent polygon $P_1'$. The set of lattice points inside $P_1'$ disjoint from $E'$ is equal to the set of lattice points inside $P_1$, and the set of lattice points inside $P_2'$ disjoint from $E'$ is equal to the set of lattice points inside $P_2$ translated by $1$ in the positive $x$-direction.
\end{proof}

To prove the intersection magnitude theorem (Theorem~\ref{thm-int-mag}), we find a triangle in each ramp-like piece that covers at least half the area, then intersect these triangles to get a smaller triangle that covers at least $\frac{1}{4}$ of the area of the smaller ramp-like piece.  Once we have this triangle in common, we need to show that it has sufficiently many lattice points.  The following lemma estimates the number of lattice points in such a triangle.

\begin{lemma}[Triangle trimming]\label{lem-trim}
Let $P$ be a triangle with at least one side parallel to an axis.  Then the number of lattice points strictly inside $P$ is at least
$$\mathrm{Area}(P) - 2\cdot\mathrm{Perimeter}(P) +1,$$
if this quantity is at least $1$.
\end{lemma}

\begin{proof}
The strategy is to use Pick's theorem, which relates the area of a lattice triangle to the number of enclosed lattice points and the number of boundary lattice points.  Our triangle $P$ does not necessarily have vertices at lattice points, so our goal is to find a large enough lattice triangle inside $P$.  First we construct a parallel line $2$ units inward from each side of $P$. We call the similar triangle defined by these parallel lines the ``trimmed triangle'', denoted $P_t$.  We show below that $P_t$ has area at least $\mathrm{Area}(P) - 2\cdot \mathrm{Perimeter}(P)$.  (In the case that there is no triangle left after the trimming process, we show that $\mathrm{Area}(P) - 2\cdot \mathrm{Perimeter}(P) < 0$ and so the lemma is vacuously true.)  If we can find an ``intermediate triangle'' $P_i$ that is a lattice triangle and is strictly between $P_t$ and $P$, then Pick's theorem states
\[\mathrm{Area}(P_i) = \#(\mathrm{interior\ lattice\ points}) + \frac{1}{2}\cdot \#(\mathrm{boundary\ lattice\ points}) - 1,\]
so we have
\[\#(\mathrm{total\ lattice\ points\ of\ }P_i) \geq \mathrm{Area}(P_i) + 1,\]
giving our goal inequality
\[\#(\mathrm{interior\ lattice\ points\ of\ }P) \geq \mathrm{Area}(P_t) + 1 \geq \mathrm{Area}(P) - 2\cdot \mathrm{Perimeter}(P) + 1.\]

We begin by estimating the area of $P_t$.  The region inside $P$ and outside $P_t$ consists of three trapezoids, each with height $2$ and one base a side of $P$.  Because the two angles bordering that side add up to less than $180^\circ$, the other base of each trapezoid---that is, the corresponding side of $P_t$---must be shorter.  Thus, the total area of the trapezoids is less than $2\cdot \mathrm{Perimeter}(P)$, giving the estimate
\[\mathrm{Area}(P_t) \geq \mathrm{Area}(P) - 2\cdot \mathrm{Perimeter}(P).\]
Suppose there is no triangle left after the trimming.  Then the inradius $r$ of $P$ is at most $2$, and connecting the vertices of $P$ to the incenter divides $P$ into three triangles, each with height $r$ and base equal to one side of $P$.  Thus we have
\[\mathrm{Area}(P) = \frac{r}{2} \cdot \mathrm{Perimeter}(P) < 2 \cdot \mathrm{Perimeter}(P),\]
and so the quantity $\mathrm{Area}(P) - 2\cdot \mathrm{Perimeter}(P)$ is negative.

\begin{figure}[h]
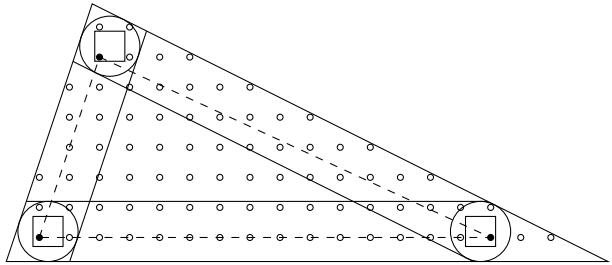

    \centering
    \figtrim
    \caption{The parallelogram trimmed off near each vertex of the original triangle contains a lattice point, because it contains a circle of radius $1$, which contains a square of side length $1$.  The lattice points from the three parallelograms form the intermediate triangle $P_i$ (dashed).}
    \label{fig:trim}
\end{figure}

At each corner of the triangle, there is a parallelogram enclosed by the two sides of the triangle and the lines parallel to each side at distance $2$, as in Figure~\ref{fig:trim}.  It suffices to find a lattice point inside each of these corner parallelograms, because these three points determine a triangle for which none of the sides crosses either a side of the original triangle $P$ or a side of the trimmed triangle $P_t$; thus, we can choose that triangle to be our intermediate triangle $P_i$.

To find the lattice point, first we observe that there is an inscribed circle of radius $1$ inside each corner parallelogram; this is because the parallelogram is the intersection of two infinite strips of width $2$, and the center lines of the two strips intersect at the center of the circle.  We also know that every square with sides parallel to the axes and of length $1$ must contain a lattice point, because tiling the plane with such squares gives lattice points at the same relative locations in each square.  Any circle of radius $1$ contains such a square of side length $1$---in fact, it contains a square of side length $\sqrt{2}$, because the diagonal has the same length $2$ as the diameter of the circle.  

Thus every corner parallelogram does contain a lattice point in its interior, so we can select one such lattice point from each corner parallelogram to define the intermediate triangle $P_i$.  Because the area of $P_i$ is greater than that of the trimmed triangle $P_t$, Pick's theorem implies that $P_i$ must have enough lattice points in its interior and boundary.
\end{proof}

Using these lemmas, we can finish proving that our pair of ramp-like polygons contains a pair of subgraphs, one on each side of the spine, that are mirror images.

\begin{proof}[Proof of intersection magnitude theorem (Theorem~\ref{thm-int-mag})]
If the common vertical edge $E$ does not have an integer $x$-coordinate, we apply the spine alignment lemma (Lemma~\ref{lem-spine}) to replace $P_1$ and $P_2$ by polygons that cut out the same subgraphs in their interiors.  Thus, we may assume that $E$ has an integer $x$-coordinate.

Let $a$ be the length of $E$, and let $b$ and $c$ be the widths of graphs $P_1$ and $P_2$, respectively.  Without loss of generality we assume $b \leq c$.  We observe that $P_1$ and $P_2$ are contained in their bounding boxes, which have $(a+1)(b+1)$ lattice points and $(a+1)(c+1)$ lattice points, respectively.  Thus it suffices to construct subgraphs $G_1$ and $G_2$ with at least $\frac{1}{20}(a+1)(b+1)$ vertices each.

Let $T_1$ and $T_2$ be right triangles constructed from the endpoints of $E$ and the vertices of $P_1$ and $P_2$ (respectively) with the greatest horizontal distance from $E$.  We construct $G_1$ by reflecting $T_2$ over $E$ and taking all the vertices in the interior of $T_1$ that are also in the interior of the reflected $T_2$; similarly, we construct $G_2$ by reflecting $T_1$ over $E$ and taking all the vertices in the interior of $T_2$ that are also in the interior of the reflected $T_1$.  Abusing notation, we let $T_1 \cap T_2$ denote the triangle formed by intersecting $T_1$ with the reflection of $T_2$.  Once we estimate its area and perimeter, we can use the triangle trimming lemma (Lemma~\ref{lem-trim}) to get a lower bound on the number of vertices of $G_1$, and hence of $G_2$ as well.

First we claim
\[\mathrm{Area}(T_1\cap T_2)\geq\frac{ab}{4}.\]
To prove this area bound, we observe that given lengths $a$, $b$, and $c$, the case where $\mathrm{Area}(T_1 \cap T_2)$ is the least is the case where the third vertex of $T_1$---that is, the vertex not on the common spine $E$---shares a $y$-coordinate with the top vertex of $E$, and the third vertex of $T_2$ shares a $y$-coordinate with the bottom vertex of $E$, or vice versa.  In this case, if $b = c$ then the width of the intersection triangle is exactly $\frac{b}{2}$ so we have $\mathrm{Area}(T_1 \cap T_2) = \frac{ab}{4}$.  If $c > b$ then the intersection triangle is larger.  Thus, in every case we have the desired area bound.

We also claim
\[\mathrm{Perimeter}(T_1\cap T_2)\leq 2(a+b).\]
This is because the perimeter of $T_1 \cap T_2$ is less than the perimeter of its bounding box, which has height $a$ and width at most $b$.

We put together the area and perimeter bounds with the triangle trimming lemma (Lemma~\ref{lem-trim}) to estimate the number of lattice points enclosed by $T_1 \cap T_2$.  It is algebraically true that for all $a,b\geq 41$, we have
$$\frac{ab}{4}-4(a+b) + 1\geq \frac{1}{20}(a+1)(b+1).$$
(To check this, we can use $ab = \frac{b}{2}a + \frac{a}{2}b > 20(a+b)$.)  Thus, using our hypothesis that $a, b, c \geq 41$, we see that the number of vertices in the interior of $T_1 \cap T_2$ satisfies the inequalities
\[\lvert G_1 \rvert \geq \frac{ab}{4}-4(a+b)+1\geq \frac{1}{20}(a+1)(b+1),\]
as desired.
\end{proof}

Having proved this estimate on the size of the mirror-image subgraphs, we can finish proving a bound on the routing number of the union of two ramp-like polygons along a shared spine.

\begin{theorem}[Routing between ramp-like]\label{thm-two-ramps}
Consider two ramp-like pieces $P_1$ and $P_2$ with common vertical spine $E$. There exists a constant $C>0$ such that $\rt(P_1\cup P_2)\leq C\cdot (w(P_1\cup P_2)+h(P_1\cup P_2)).$
\end{theorem}

\begin{proof}
First consider the case where $w(P_1)$, $w(P_2)$, and the length of $E$ are all at least $41$, so the intersection magnitude theorem (Theorem~\ref{thm-int-mag}) applies.  Fix a home configuration of tokens on the vertices of $P_1 \cup P_2$, and consider an arbitrary starting configuration.  The ramp-like bound (Theorem~\ref{thm-ramplike}) implies that we can efficiently route the tokens within $P_1$ and the tokens within $P_2$.  Thus, what we need to show is that we can efficiently route the tokens into their home halves---that is, those that belong in $P_1$ should go to $P_1$ and those that belong in $P_2$ should go to $P_2$.

We label each token either black or white indicating whether it belongs in $P_1$ or $P_2$, respectively, in the home configuration.  If the common spine $E$ has an integer $x$-coordinate, then some lattice points are shared between $P_1$ and $P_2$.  In this case we count those lattice points as part of $P_1$ and not $P_2$, so that each token belongs in exactly one of the halves, and without loss of generality, we may assume that there are at least as many lattice points in $P_1$ as in $P_2$.  In any configuration, we say that a given token is \textbf{\textit{improper}} if it is in the opposite half from where it belongs.  The number of improper tokens in $P_1$ is always equal to the number of improper tokens in $P_2$, which is at most the total number of lattice points in $P_2$.

We use the intersection magnitude theorem (Theorem~\ref{thm-int-mag}) to find subgraphs $G_1$ and $G_2$ in $P_1$ and $P_2$ that are reflections over a vertical line and have size at least $\frac{1}{20}\lvert P_2 \rvert$.  Then we can move up to $\lvert G_1\rvert = \lvert G_2 \rvert$ improper tokens into their home halves, using the following sequence of phases:
\begin{enumerate}
\item Use the ramp-like bound (Theorem~\ref{thm-ramplike}) to route within $P_1$ and $P_2$ separately so that as many improper tokens as possible are in $G_1$ and $G_2$.  If there are at least $\lvert G_1\rvert$ improper tokens on each side, then $G_1$ and $G_2$ become completely filled with improper tokens.
\item In the case where $G_1$ and $G_2$ do not become completely filled with improper tokens, continue to route within $P_2$ so that the locations of the improper tokens in $G_2$ are exactly the mirror image of the locations of the improper tokens in $G_1$.
\item Route each row of $P_1 \cup P_2$ simultaneously so that the improper tokens in $G_1$ exchange places with their mirror-image improper tokens in $G_2$, leaving no more improper tokens in either $G_1$ or $G_2$.
\end{enumerate}

Each of these three phases takes at most $C(w(P_1 \cup P_2) + h(P_1 \cup P_2))$ routing steps, for some constant $C$.  Repeating up to $20$ times if necessary, we can move every token into its home half so that no improper tokens remain.  Applying the ramp-like bound (Theorem~\ref{thm-ramplike}) once more to route within each half, we move all tokens to their home lattice points in at most $C(w(P_1 \cup P_2) + h(P_1 \cup P_2))$ routing steps, for some constant $C$.

We now consider the case where the widths of $P_1$ and $P_2$ and the length of $E$ are not all at least $41$. Suppose without loss of generality that it is $P_2$ that has height or width less than $41$.  We apply Lemma~\ref{lem-hair} with $G = P_1 \cup P_2$, $K = P_1$, the constant $c_1$ is $41$, and $S$ is the rightmost column of $P_1$.  Because we have a bound on the routing number of $P_1$, Lemma~\ref{lem-hair} states that the routing number of $P_1 \cup P_2$ is at most $C(w(P_1 \cup P_2) + h(P_1 \cup P_2))$ for some constant $C$.
\end{proof}

Using the bound for a pair of ramp-like polygons, we can follow a similar sequence of steps again to finish proving the bound on routing number of burger bun polygons.

\begin{theorem}[Burger bun bound]\label{thm-burger-bun}
There exists a constant $C>0$ such that for any burger bun polygon $P$, the routing number of $P$ satisfies the bound $\rt(P)\leq C\cdot (w(P)+h(P))$.
\end{theorem}

\begin{proof}
Let $E$ be the spine of $P$.  Without loss of generality we may assume that $E$ is vertical, so $E$ divides $P$ into a left half $P_1$ and a right half $P_2$.  Each of $P_1$ and $P_2$, if it is not ramp-like already, is the union of two ramp-like pieces sharing a horizontal spine.  Thus, Theorem~\ref{thm-two-ramps} gives a bound on the routing number of $P_1$ and $P_2$ separately.  In the present proof, we follow the proof of Theorem~\ref{thm-two-ramps}, but instead of using the ramp-like bound (Theorem~\ref{thm-ramplike}) to route the two ramp-like halves, we use Theorem~\ref{thm-two-ramps} itself to route $P_1$ and $P_2$.

We still need to prove an analogue of the intersection magnitude theorem (Theorem~\ref{thm-int-mag}) that applies to the present $P_1$ and $P_2$, which are not necessarily ramp-like.  To do this, we construct two right triangles $T_1$ and $T_2$ that are ramp-like with common spine $E$, such that the intersection of $T_1$ with the reflection of $T_2$ is contained in the intersection of $P_1$ with the reflection of $P_2$.  Then we apply the intersection magnitude theorem (Theorem~\ref{thm-int-mag}) to $T_1$ and $T_2$.

\begin{figure}[h]
    \centering
    \begin{tikzpicture}[scale=.7]

\fill[gray!3!] (0,0)--(0,6)--(-3,5)--(-2,3)--(0,0);
\draw[stroke=thick] (0,0)--(0,6)--(-1.5,4.8)--(-2.5,2.1)--(0,0);
\fill[gray!9!](0,0)--(0,6)--(-1.5,4.8)--(-2.5,2.1)--(0,0);
\draw[stroke=thick] (0,0)--(0,6)--(3,5)--(2,2.5)--(0,0);
\fill[gray!9!] (0,0)--(0,6)--(3,5)--(2,2.5)--(0,0);
\draw[stroke=thick] (0,6)--(-3,5)--(-2,2.5)--(0,0);

\draw[dashed] (0,6)--(-2.5, 2.1) (0,0)--(3.846, 0)--(0, 6);
\draw[dashed] (0,6)--(-3.6,6)--(0,0);

\filldraw (0,0) circle (2pt);
\filldraw (0,6) circle (2pt);
\filldraw (-3.6,6) circle (2pt);
\filldraw (3.846, 0) circle (2pt);
\filldraw (-3,5) circle (2pt);
\filldraw (-2.5,2.1) circle (2pt);
\filldraw (-1.8595,3.09917) circle (2pt);

\node at (-.7,3.2) {$P_1$};
\node at (.7, 3.2) {$P_2$};
\node at (0,-.3) {$B$};
\node at (0,6.3) {$A$};
\node at (-3.6,6.3) {$W_1$};
\node at (3.846,-.3) {$W_2$};
\node at (3.4,5) {$V_2$};
\node at (-2.8,2.1) {$V_1$};
\node at (-1.5595,3.09917) {$C$};

\end{tikzpicture}
    \caption{We use the intersection magnitude bound for ramp-like polygons to show a similar bound for burger bun halves, by applying it to ramp-like triangles $T_1 = ABW_1$ and $T_2 = ABW_2$, which are sufficiently large but still have (reflected) intersection inside our burger bun halves.}
    \label{fig:burgers}
\end{figure}
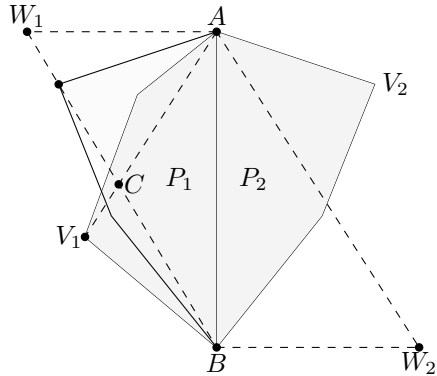

We construct $T_1$ and $T_2$ as follows, shown in Figure~\ref{fig:burgers}. Let $V_1$ and $V_2$ be the points on $P_1$ and $P_2$ (respectively) farthest from edge $E$, and let $A$ and $B$ be the top and bottom points of $E$.  Then the triangle $ABV_1$ intersects the reflection of $ABV_2$ in some triangle $ABC$.  We construct the third vertex $W_1$ of $T_1$ (that is, the vertex other than $A$ and $B$) by continuing the segment $BC$ upward until it intersects the horizontal line containing $A$.  Similarly, we construct the third vertex $W_2$ of $T_2$ by continuing the reflection of $AC$ downward until it intersects the horizontal line containing $B$.

By applying the intersection magnitude theorem (Theorem~\ref{thm-int-mag}) to $T_1$ and $T_2$, we find subgraphs $G_1$ and $G_2$ in $T_1$ and $T_2$, consisting of at least $\frac{1}{20}\min\{\lvert T_1\rvert, \lvert T_2\rvert\}$ vertices each.  The construction of $T_1$ and $T_2$ guarantees that $G_1$ is also in $P_1$ and that $G_2$ is also in $P_2$.  We compare $\lvert T_1 \rvert$ and $\lvert T_2\rvert$ to $\lvert P_1\rvert$ and $\lvert P_2 \rvert$ by observing that $T_1$ and $T_2$ each have at least half as many vertices as their bounding boxes, and that those bounding boxes contain the bounding boxes of $P_1$ and $P_2$.  Thus we have $\lvert T_1 \rvert \geq \frac{1}{2}\lvert P_1\rvert$ and $\lvert T_2 \rvert \geq \frac{1}{2}\lvert P_2 \rvert$.  The result is the lower bound
\[\lvert G_1\rvert = \lvert G_2 \rvert \geq \frac{1}{40} \min\{\lvert P_1 \rvert, \lvert P_2 \rvert\},\]
in the case where the width and height of $P_1$ and $P_2$ are all at least $41$.

Thus, following the proof of the routing between ramp-like theorem (Theorem~\ref{thm-two-ramps}), in the case where the height and width of $P_1$ and $P_2$ are all at least $41$, we can move at least $\frac{1}{40}$ of the improper tokens into their home halves by first routing each half to put the improper tokens into $G_1$ and $G_2$, and then routing each row simultaneously to swap the improper tokens in $G_1$ with the improper tokens in $G_2$.  Repeating this process at most $40$ times puts every token into its home half, and then one more instance of routing within the two halves moves every token to its home vertex.

In the case where the height or width of (without loss of generality) $P_2$ is less than $41$, as in the proof of Theorem~\ref{thm-two-ramps} we may apply Lemma~\ref{lem-hair} with $G = P_1 \cup P_2$, $K = P_1$, $c_1 = 41$, and $S$ is the rightmost column of $P_1$.
\end{proof}

\section{Proof of main theorem}\label{sec-main-proof}

In this section we prove the main theorem.  For convenience we reproduce the statement of the main theorem here.

\begin{reptheorem}{thm-main}
Let $P$ be a connected convex grid piece. Then the routing number of $P$ satisfies the bound $\rt(P) \leq  C (w(P) + h(P))$ for some universal constant $C$. 
\end{reptheorem}

The idea of the proof is much simpler than the details.  We show that our arbitrary convex polygon $P$ is related to a burger bun polygon by a shear transformation of at most $45^\circ$.  Roughly, this transformation corresponds to a map between the sets of enclosed lattice points that stretches distances by at most a fixed factor.  We show that if two graphs are related by a map that stretches by at most a fixed factor, then their routing numbers are also related by at most a fixed factor.  Thus, if we shear the original polygon to get a burger bun, then routing the resulting burger bun polygon helps us to route the original polygon.

\begin{figure}[h]
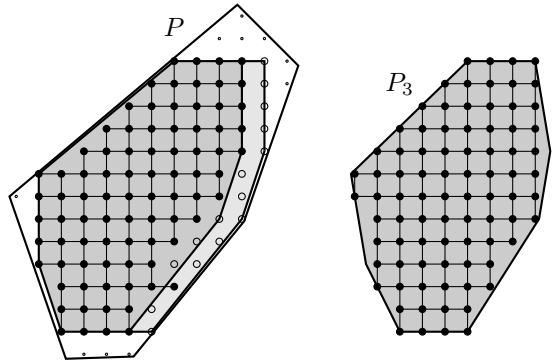

    \centering
    \figsummary
    \caption{We cut off any short rows or columns of the original polygon $P$ to form $P_1$ (left, shaded), then cut off the rightmost vertex in each row to form $P_2$ (left, darkly shaded) before shearing $P_2$ to get a burger bun polygon $P_3$.  The graph $\psi(P_3)$ inside $P_1$ has the same number of vertices per row as $P_3$ has.}
    \label{fig:main-summary}
\end{figure}

The actual proof becomes more complicated to account for how the shear transformation does not respect the integer lattice---in particular, it does not necessarily preserve the number of lattice points inside the polygon.  Lemmas~\ref{lem-trimp} and~\ref{lem-overlap} describe how we cut off part of $P$ to form $P_1 \subseteq P$, and Lemma~\ref{lem-p2} describes how we cut off a little more to form $P_2 \subseteq P_1$.  Lemmas~\ref{lem-shear} and~\ref{lem-bb-conn} describe how we shear $P_2$ to get a burger bun polygon $P_3$, which we know how to route.  Lemma~\ref{lem-psi-lip} describes how to map the lattice points inside $P_3$ into $P_1$ to form a subgraph $\psi(P_3)$ of $P_1$, and the bounded stretch theorem (Theorem~\ref{thm-lip}) implies that the routing number of $\psi(P_3)$ is at most a constant factor greater than that of $P_3$.  Then Lemma~\ref{lem-psi-bdry} checks the hypotheses of Lemmas~\ref{lem-hair} and~\ref{lem-skin}, which will show that the routing number of $P$ is not much greater than that of $\psi(P_3)$.  Figure~\ref{fig:main-summary} shows the relationship between these polygons and their associated graphs.

Before shearing, we want to know that the shear does not affect whether the enclosed lattice graph is connected.  To do this, in the next two lemmas we trim off the short rows and columns of $P$ that would be at risk of being pulled apart by the shear.

\begin{lemma}[Constructing $P_1$]\label{lem-trimp}
Let $P$ be a convex polygon enclosing strictly more than $4(w(P) + h(P))$ lattice points.  Then there exists a polygon $P_1 \subseteq P$ such that
\begin{itemize}
\item $P_1$ is the convex hull of its enclosed lattice points;
\item The subgraph $P \setminus P_1$ has at most $4(w(P) + h(P))$ vertices; and
\item Every row and column of $P_1$ has at least $4$ lattice points.
\end{itemize}
\end{lemma}

\begin{proof}
We consider the top and bottom rows of $P$ and the leftmost and rightmost columns of $P$.  If each of these has at least $5$ lattice points, we set $P_1$ to be the convex hull of the lattice points in $P$.  Otherwise, we iteratively remove one row or column at a time from the graph, choosing either the top row, the bottom row, the leftmost column, or the rightmost column, whichever has at most $4$ lattice points.  Once the top and bottom rows and the leftmost and rightmost columns of the remaining graph all have at least $5$ lattice points, we set $P_1$ to be the convex hull of all the lattice points remaining.

Each deletion step reduces the number of lattice points by at most $4$, while also reducing either the width or the height (possibly both) by $1$.  Thus $P \setminus P_1$ has at most $4(w(P) + h(P))$ lattice points, and in particular $P_1$ is nonempty.

To show that every row and column has at least $4$ lattice points, because of the symmetry it suffices to show that every row has at least $4$ lattice points.  Consider the parallelogram formed by any choice of $5$ consecutive lattice points from the top row of $P_1$ and $5$ consecutive lattice points from the bottom row of $P_1$.  Because $P_1$ is convex, this parallelogram is contained in $P_1$.  Every horizontal cross-section of the parallelogram is an interval of length $4$, so each cross-section at integer height must contain either $4$ or $5$ lattice points.  Thus every row has at least $4$ lattice points.
\end{proof}

The following lemma shows that consecutive rows of $P_1$ are well connected to each other, as are consecutive columns.

\begin{lemma}[Property of $P_1$]\label{lem-overlap}
Let $P_1$ be a convex polygon with at least $4$ lattice points in every row and column.  Then every two consecutive rows of $P_1$ have at least $3$ columns in common, and every two consecutive columns of $P_1$ have at least $3$ rows in common.
\end{lemma}

\begin{proof}
Because of the symmetry, it suffices to show that every two consecutive rows of $P_1$ have at least $3$ columns in common.  Suppose for the sake of contradiction that there are two consecutive rows with at most two columns in common; without loss of generality, suppose that it is the two (or more) rightmost vertices in the upper row that are not adjacent to vertices in the lower row.  Let $u$ be the vertex in the upper row just to the right of the shared columns, and let $v$ be the vertex in the lower row just to the left of the shared columns, as in Figure~\ref{fig:overlap}.

\begin{figure}[h]
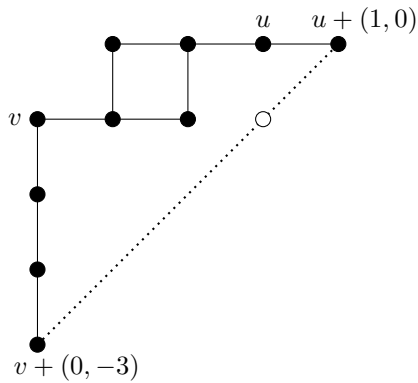

    \centering
    \figoverlap
    \caption{If the rows of $u$ and $v$ have fewer than $3$ columns in common, we can contradict the convexity of $P_1$.}
    \label{fig:overlap}
\end{figure}

Then $v$ has no vertex immediately above it, so it must have at least three vertices below it in the same column, which means that the lattice point $v + (0, -3)$ is a vertex in $P_1$.  We also know that $u + (1, 0)$ is a vertex in $P_1$.  The segment with endpoints $v + (0, -3)$ and $u + (1, 0)$ has slope at least $1$, because the $x$-coordinates of $u$ and $v$ differ by at most $3$.  The lattice point $u + (0, -1)$ lies on or to the left of this segment, but we have assumed that it is not in $P_1$, giving a contradiction.  Thus we may conclude that each pair of consecutive rows has more than two columns in common.
\end{proof}

It turns out that we want to trim off the right side of $P_1$ to form $P_2$ before shearing.  This ensures that later when we construct $\psi(P_3)$, it fits inside $P_1$.

\begin{lemma}[Constructing $P_2$]\label{lem-p2}
Let $P_1$ be a convex polygon.  Then there is a convex polygon $P_2$ that encloses all of the lattice points inside $P_1$ except the rightmost lattice point of each row.
\end{lemma}

\begin{proof}
Let ${P_1}_{\mathrm{ymax}}$ and ${P_1}_{\mathrm{ymin}}$ be points on $P_1$ with the maximum and minimum $y$-coordinate, respectively. Consider the portion of the boundary of $P_1$  that is between the points ${P_1}_{\mathrm{ymax}}$ and ${P_1}_{\mathrm{ymin}}$ when moving clockwise from ${P_1}_{\mathrm{ymax}}$ to ${P_1}_{\mathrm{ymin}}$. Translate this piecewise-linear curve to the left by $1$, and define the new polygon $P_2$ to be the subset of $P_1$ that is to the left of this
translated side. Then $P_2$ has exactly one less vertex per row than $P_1$.
\end{proof}

We construct a shear that transforms $P_2$ into a burger bun polygon, and label this burger bun polygon $P_3$, in the following lemma.

\begin{lemma}[Constructing $P_3$]\label{lem-shear}
Let $P_2$ be a convex polygon with $w(P_2) \leq h(P_2)$. Then there exists a shear $$S=\begin{bmatrix} 1 & m \\ 0 & 1\end{bmatrix}$$ with $|m| \leq 1$ such that the resulting polygon $P_3 = SP_2$ is burger bun, and we have $w(P_3) \leq 2w(P_2)$.
\end{lemma}\label{shearing}

\begin{proof}

Let $p = P_{\mathrm{ymax}}$ and $q = P_{\mathrm{ymin}}$. Then $\big\lvert\frac{q_x-p_x}{p_y-q_y}\big\rvert \le 1$ because $|q_x-p_x| \le w(P) \le h(P)=|p_y-q_y|$.
 Define the horizontal shear
\[ S = \begin{bmatrix}
1 &\frac{q_x-p_x}{p_y-q_y}\\
0 & 1
\end{bmatrix}.\]

The region $SP= P'$ is a convex polygon since $S$ is linear. Moreover, $P'$ is burger bun because $S$ fixes the $y$-coordinate of each point and we have constructed the matrix $S$ so that $x(P'_{\mathrm{ymax}}) = x(P'_{\mathrm{ymin}})$, as follows:
\begin{align*}
x(P'_{\mathrm{ymax}}) = x(Sp) &= p_x+p_y\frac{q_x-p_x}{p_y-q_y}
=\frac{-p_xq_y + p_yq_x}{p_y-q_y}=\\
&=q_x+q_y\frac{q_x-p_x}{p_y-q_y} = x(Sq) = x(P'_{\mathrm{ymin}}).
\end{align*}

To show that $w(P_3) \leq 2w(P_2)$, let $a$ and $b$ be two arbitrary points in $P_2$.  Then the $x$-coordinates of their images in $P_3$ are $a_x + ma_y$ and $b_x + mb_y$, which have absolute difference at most 
\[\abs{a_x - b_x} + \abs{m}\cdot \abs{a_y - b_y} \leq w(P_2) + \frac{w(P_2)}{h(P_2)} \cdot h(P_2) = 2w(P_2).\]
\end{proof}

Before we can apply our burger-bun routing theorem to $P_3$, we need to check in the following lemma that $P_3$ is connected.

\begin{lemma}[Property of $P_3$]\label{lem-bb-conn}
Let $P_3$ be a burger bun polygon with vertical spine. Suppose that $P_3$ contains at least $3 \cdot w(P_3)$ vertices in total and at least $2$ vertices in each row. Then the graph $P_3$ is connected.
\end{lemma}

\begin{proof}
If the spine is at an integer $x$-coordinate, then $P_3$ is connected because every vertex is connected by a horizontal path to the spine.  Otherwise, we consider the two columns surrounding the spine.  Every vertex in $P_3$ can be connected by a horizontal path to one of these two columns, so it suffices to show that the two columns have a row in common.  Suppose to the contrary that they do not; without loss of generality, the $y$-coordinates of the subgraph of $P_3$ to the left of the spine are all greater than the $y$-coordinates of the subgraph to the right.

\begin{figure}[h]
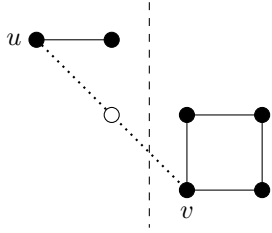

    \centering
    \figthrows
    \caption{If a burger bun polygon has at least $3$ rows, with at least $2$ vertices in each row, then the subgraphs on either side of the spine must connect to each other, or else contradict the convexity of the polygon.}
    \label{fig:3rows}
\end{figure}

Because $P_3$ has at least $3\cdot w(P_3)$ vertices, it must have at least three rows; without loss of generality, the left subgraph has at least one row, and the right subgraph has at least two rows.  Consider the second-to-right vertex $u$ in the bottom row of the left subgraph, and the second-to-top vertex $v$ in the leftmost column of the right subgraph, as in Figure~\ref{fig:3rows}.  The midpoint of $u$ and $v$ would connect the left subgraph to the right subgraph, so our assumptions imply that it is not a vertex in $P_3$; however, this contradicts the convexity of $P_3$.  Thus, it is impossible for the two columns surrounding the spine not to have a row in common, and so $P_3$ is connected.
\end{proof}

The burger bun bound (Theorem~\ref{thm-burger-bun}) shows that we can route $P_3$.  To show that this helps us route $P$, we start by finding a distorted copy of $P_3$ inside $P_1$.

\begin{lemma}[Constructing $\psi(P_3)$]\label{lem-psi-lip}
Let $P_1$, $P_2$ and $P_3$ be convex polygons constructed in Lemmas~\ref{lem-trimp}, \ref{lem-p2}, and \ref{lem-shear}.  Then there is an injective map $\psi$ from the vertex set of $P_3$ to the vertex set of $P_1$ that sends adjacent vertices in $P_3$ to vertices no more than 3 edges away from each other in $P_1$.  In particular, its image $\psi(P_3)$ is connected.
\end{lemma}

\begin{proof}
First we construct the map $\psi$.  We start by defining new coordinates for the vertices in $P_1$ and $P_3$.  Without loss of generality suppose that the least $y$ coordinate of vertices in each polygon is $1$. Let \(r_i\) be the set of vertices in \(P_1\) with \(y\) coordinate \(i\), so $r_i$ is the $i$th row of $P_1$. Let \(r_{i,j}\) be the $j$th vertex in row \(r_i\), counting from left to right. That is, $r_{i,1}$ is the vertex with least $x$ coordinate in row $i$, and $r_{2, i}$ has the second least $x$ coordinate and so on.   Similarly, label the vertices of $P_3$ as $r'_{i,j}$.  Define an injection \(\psi\)  from the vertex set of $P_3$ to $P_1$, given by \(\psi(r'_{i,j}) = r_{i,j}\). That is, $\psi$ sends the $j$th vertex in row $i$ of $P_3$ to the $j$th vertex in row $i$ of $P_1$. We have constructed the polygons so that each row of $P_3$ has no more vertices than the corresponding row of $P_1$: every horizontal cross-section of $P_3$ has the same length as the corresponding cross-section of $P_2$, so the number of lattice points in each row of $P_3$ differs by at most $1$ from the number of lattice points in the corresponding row of $P_2$, and we know that $P_1$ has one more lattice point per row than $P_2$ has.  Thus, our map $\psi$ is well-defined.

Before proving that $\psi$ stretches distances by at most a factor of $3$, as a preliminary step we show that every two consecutive rows of $\psi(P_3)$ have at least one column in common.  By Lemma~\ref{lem-overlap} we know that every two consecutive rows of $P_1$ have at least three columns in common.  Then $P_2$ is like $P_1$, but with the rightmost vertex deleted from each row, so every two consecutive rows of $P_2$ have at least two columns in common.  Then because corresponding rows of $P_2$ and $P_3$ differ in length by at most $1$, we know that $\psi(P_3)$ is like $P_1$, but with up to two of the rightmost vertices deleted from each row, so every two consecutive rows of $\psi(P_3)$ have at least one column in common.

Using these common columns of pairs of consecutive rows, we can prove that $\psi$ sends adjacent vertices to vertices no more than $3$ edges away from each other.  Adjacent vertices in a row of $P_3$ are sent to adjacent vertices, so we only need to check what $\psi$ does to adjacent vertices in a column of $P_3$.  Thus, it suffices to show that each pair of adjacent rows do not shift in relation to each other by more than $2$ edges under $\psi$. That is, for all rows $r_j$ and $r_{j+1}$ in $P_3$, we must show that the difference $|[x(r_{j,1})- x(r_{j+1,\text{ }1})] - [x(r'_{j,1})-x(r'_{j+1,\text{ }1})]|$ is no more than 2. 

Let $(x_1,j)$ be the leftmost point in $P_3$ with $y$-coordinate $j$, and $(x_2,j+1)$ be the leftmost point of $P_3$ with $y$-coordinate $j+1$. Note that $r_{j,1}$ and $r_{j+1,\text{ }1}$ are the first vertices in each of these rows, so $r_{j,1} = (\lceil{x_1}\rceil,j)$ and $r_{j+1,\text{ }1} = (\lceil{x_2}\rceil,j+1)$. After the shear by the matrix 
\begin{equation*}
M = \begin{bmatrix} 1 & m \\ 0 & 1  \end{bmatrix},
\end{equation*}
the points $(x_1,j)$ and $(x_2,j+1)$ are sent to $(x_1 + mj, j)$ and $(x_2 + m(j+1), j+1)$. Therefore, the first vertices in these rows of $P_1$, namely $r'_{j,1}$ and $r'_{j+1,\text{ }1}$, have coordinates $(\lceil{x_1 + mj}\rceil, j)$ and $(\lceil{x_2+m(j+1)}\rceil, j+1)$ respectively. Bounding the distance between these two points in $P_1$, we have
\begin{equation}
\begin{aligned}
 (x_1 - x_2) - |m| - 1 <  \lceil{x_1 + mj}\rceil - \lceil{x_2+m(j+1)}\rceil <  (x_1 - x_2)  + |m| + 1.
\end{aligned}
\end{equation}
And, bounding the distance between the original points in $P_3$, we have
\begin{equation}
(x_1-x_2) - 1 < \lceil{x_1}\rceil - \lceil{x_2}\rceil < (x_1 - x_2) + 1.
\end{equation}

Combining these inequalities, we find that the absolute difference between the quantities $\left\lceil x_1 \right\rceil - \left\lceil x_2 \right\rceil$ and $\left\lceil x_1 + mj \right\rceil - \left\lceil x_2 + m(j+1)\right\rceil$ is strictly less than $\abs{m} + 2$, and thus is strictly less than $3$.  Because the difference is an integer, it must be at most $2$; in other words, the rows shift by no more than $2$ vertices away from each other. Thus, the map $\psi$ sends adjacent vertices to vertices no more than $3$ edges away from each other. 
\end{proof}

Applying the next theorem to $\psi$ gives us a way to route $\psi(P_3)$ using our knowledge of how to route $P_3$.

\begin{theorem}[Bounded stretch]\label{thm-lip}
Let $A$ and $B$ be lattice graphs in $\mathbb{Z} \times \mathbb{Z}$, and let $\psi$ be a bijection between the vertices of $A$ and the vertices of $B$, such that for any pair of adjacent vertices $v_1, v_2 \in A$, the path-length distance between $\psi(v_1)$ and $\psi(v_2)$ in $B$ is at most a constant $c$. Then $\rt(B) \le C(c) \cdot \rt(A)$, where $C(c)$ is a constant depending only on $c$.
\end{theorem}

\begin{proof}
We would like to color the edges of $A$, such that if $(v_1, v_2)$ and $(v'_1, v'_2)$ are two edges of the same color, and we draw shortest paths between $w_1 = \psi(v_1)$ and $w_2 = \psi(v_2)$ and between $w'_1 = \psi(v'_1)$ and $w'_2 = \psi(v'_2)$ in $B$, then these two paths are disjoint.  To do this, we would like to assign a color to each pair of vertices in $B$ that are within distance $c$ of each other, such that pairs of the same color are more than distance $c$ apart.  We do this by first coloring the vertices of $B$, and then coloring the distance $c$ pairs by the color pairs of their vertices.

Our first coloring assigns colors to the vertices of $B$, such that if two vertices have the same color, they have distance greater than $2c$ in $B$.  To do this, we construct a graph $B'$ with the same vertex set as $B$, with an edge between vertices $w_1$ and $w_2$ whenever their distance is at most $2c$. The maximal degree of any vertex in $B'$ is at most $4c(2c+1)$, since there are $4i$ lattice points with distance exactly $i$ away in $\mathbb{Z} \times \mathbb{Z}$ for each $1 \leq i \leq 2c$, and if a pair of vertices have distance at most $2c$ in $B$, they also have distance at most $2c$ in $\mathbb{Z} \times \mathbb{Z}$. Therefore, using a greedy strategy we can color the vertices of $B'$ with no more than $4c(2c+1) + 1$ colors, so that no two vertices of the same color are adjacent in $B'$. 

Our second coloring has one color for each pair of colors in the first coloring, for a total of $\binom{4c(2c+1)+1}{2} = (4c(2c+1)+1)\cdot 2c(2c+1)$ colors.  It assigns one color to each pair of vertices $w_1, w_2$ in $B$ that have distance at most $c$, given by the pair of colors of $w_1$ and $w_2$ in the first coloring.  Suppose that $(w_1, w_2)$ and $(w'_1, w'_2)$ are two pairs of vertices at distance at most $c$ in $B$, and they have the same color in the second coloring; without loss of generality this means that $w_1$ and $w'_1$ have the same color in the first coloring, as do $w_2$ and $w'_2$.  

We claim that if we choose any shortest path from $w_1$ to $w_2$ in $B$, and any shortest path from $w'_1$ to $w'_2$ in $B$, then these two paths are disjoint.  Suppose to the contrary that these two paths cross at a point. Then either $w_1$ or $w_2$ is within $\frac{1}{2}c$ of the common point, as is either $w'_1$ or $w'_2$.  This implies that one of the pairs $(w_1, w'_1)$, $(w_1, w'_2)$, $(w_2, w'_1)$, or $(w_2, w'_2)$ has distance at most $c$, and so either $w_1$ and $w'_1$ have distance at most $2c$, or $w_2$ and $w'_2$ have distance at most $2c$.  This contradicts their having been colored the same color in the first coloring.  Therefore, any shortest path from $w_1$ to $w_2$ in $B$, and any shortest path from $w'_1$ to $w'_2$ in $B$ are disjoint.

For each edge in $A$, we color it by taking the corresponding pair of vertices in $B$, and finding the color of that pair in the second coloring.  Then, given a pair of configurations in $B$, we use the following method to route between them. Take the corresponding pair of configurations in $A$, and consider a shortest sequence of steps to route between them. For each step in $A$, a set of disjoint swaps, all swaps along edges of the same color in $A$ can be carried out in parallel in $B$, because the corresponding paths are disjoint. There are at most $(4c(2c+1)+1)\cdot 2c(2c+1)$ colors of edges in $A$, and each color may take $c + 1$ steps to route in $B$, since a path of length $c$ can be routed in $c + 1$ steps. Therefore, it will take at most 
\[C(c) = (4c(2c+1)+1)\cdot 2c(2c+1)\cdot (c+1)\] 
steps in total to route a single step in $A$, a set of disjoint swaps. Therefore, $\rt(B) \leq C(c) \cdot \rt(A)$.
\end{proof}

At this point, we have obtained a bound on the routing number of $\psi(P_3)$, using the bound on the routing number of $P_3$.  To extend this result to route all of $P$, we have to check the hypotheses of Lemmas~\ref{lem-hair} and~\ref{lem-skin}, and then applying these lemmas will show that because $P$ is not too much bigger than $\psi(P_3)$, its routing number is also not too much bigger.

\begin{lemma}[Hair and skin of $\psi(P_3)$]\label{lem-psi-bdry}
Let $P$ be a convex polygon enclosing strictly more than $4(w(P) + h(P))$ lattice points, let $P_1$, $P_2$, and $P_3$ be convex polygons constructed in Lemmas~\ref{lem-trimp}, \ref{lem-p2}, and \ref{lem-shear}, and let $\psi \colon\thinspace P_3 \rightarrow P_1$ be the injection constructed in Lemma~\ref{lem-psi-lip}.  Then the subgraph $\psi(P_3)$ of $P$ has the following properties:
\begin{enumerate}
\item $P \setminus \psi(P_3)$ has at most $6(w(P) + h(P))$ vertices; and
\item There is a subset $S$ of vertices in $\psi(P_3)$ such that $S$ has at most $4(w(P) + h(P))$ vertices and the induced subgraph $S \cup (P \setminus \psi(P_3))$ of $P$ is connected.
\end{enumerate}
\end{lemma}

\begin{proof}
Let $P_4$ be the convex polygon obtained from $P_2$ by removing the rightmost vertex of each row, as in Lemma~\ref{lem-p2}.  Then $\psi(P_3)$ contains $P_4$.  We know that $P_1 \setminus P_4$ contains exactly $2h(P_1)$ vertices and that $P \setminus P_1$ contains at most $4(w(P) + h(P))$ vertices, so in total $P \setminus P_4$, and therefore $P \setminus \psi(P_3)$, contains at most $4w(P) + 6h(P)$ vertices.

Lemma~\ref{lem-skin} implies that $P_4$ has at most $2(w(P_4) + h(P_4)) \leq 2(w(P) + h(P))$ vertices in its boundary.  We choose $S$ to contain the boundary of $P_4$, as well as all of $\psi(P_3) \setminus P_4$.  Because $\psi(P_3)$ has at most two more vertices in each row than $P_4$, the number of vertices in $S$ is at most $2w(P) + 4h(P)$.  

To show that $S \cup (P \setminus \psi(P_3))$ is connected, we observe that it is the same induced subgraph of $P$ as the union of the boundary of $P_4$ with $P \setminus P_4$; Lemma~\ref{lem-skin} states that because $P$ is connected and contains $P_4$, this induced subgraph is also connected.
\end{proof}

Finally we are ready to finish proving the bound on routing number of arbitrary convex polygons.

\begin{proof}[Proof of Theorem~\ref{thm-main}]
Let $P$ be a convex polygon such that the grid piece contained in $P$ is connected.  If $P$ has at most $4(w(P) + h(P))$ vertices, then the tree bound (Theorem~\ref{thm-tree}) implies that $\rt(P) \leq 12(w(P) + h(P))$, and there is nothing more to prove.  Thus, we may assume that $P$ has more than $4(w(P) + h(P))$ vertices.

We apply Lemma~\ref{lem-trimp} to find $P_1$ inside $P$ with at least $4$ vertices per row and column.  Without loss of generality, we may assume that $h(P_1) \geq w(P_1)$.  Then we apply Lemma~\ref{lem-p2} to find $P_2$ inside $P_1$ by removing the rightmost vertex of each row.  Then $P_2$ has the same height as $P_1$, and $w(P_2) = w(P_1) - 1$.  We apply Lemma~\ref{lem-shear} to shear $P_2$ to get a burger bun polygon $P_3$.  

We would like to apply Lemma~\ref{lem-bb-conn} to check that $P_3$ is connected, so we need to estimate the number of vertices in $P_3$.  We know that $P_1$ has more than $4(w(P_1) + h(P_1))$ vertices, so because $P_2$ is missing one vertex from each row, we see that $P_2$ has more than $4w(P_1) + 3h(P_1)$ vertices.  Then every horizontal cross-section of $P_3$ has the same length as the corresponding cross-section of $P_2$, so the number of vertices in each row of $P_3$ differs by at most $1$ from the number of vertices in the corresponding row of $P_2$.  This implies that $P_3$ has more than $4w(P_1) + 2h(P_1)$ vertices.  Lemma~\ref{lem-shear} tells us that $w(P_3) \leq 2w(P_2)$, and we know that $w(P_2) \leq w(P_1) \leq h(P_1)$, so the number of vertices in $P_3$ is more than $4w(P_2) + 2w(P_2) \geq 3w(P_3)$.  Thus we may apply Lemma~\ref{lem-bb-conn} to conclude that $P_3$ is connected.

The burger bun bound (Theorem~\ref{thm-burger-bun}) implies that because $P_3$ is connected and burger bun, we have $\rt(P_3) \leq C(w(P_3) + h(P_3))$ for some constant $C$.  Then Lemma~\ref{lem-psi-lip} and the bounded stretch theorem (Theorem~\ref{thm-lip}) together imply that $\rt(\psi(P_3)) \leq C(w(P_3) + h(P_3))$ for some larger constant $C$, and so $\rt(\psi(P_3)) \leq C(2w(P_1) + h(P_1))$.  Lemma~\ref{lem-psi-bdry} implies that Lemma~\ref{lem-hair} applies to $G = P$ and $K = \psi(P_3)$, so we may conclude that $\rt(P) \leq C(w(P) + h(P))$ for some constant $C$, as desired.
\end{proof}

\bibliographystyle{amsalpha}
\bibliography{mathilyest-references}{}
\end{document}